\documentclass[10pt]{article}

\usepackage{amssymb}
\usepackage{amsmath}
\usepackage{amsthm}
\usepackage{graphicx}
\usepackage{siunitx}

\newcommand{\R}{\mathbb R}
\newcommand{\C}{\mathbb C}
\newcommand{\V}{\mathbb V}
\newcommand{\St}{\mathbb S}
\newcommand{\Oh}{{\cal O}}
\newcommand{\st}{{\triangle t}}
\newcommand{\hQ}{{\hat{Q}}}
\newcommand{\hR}{{\hat{R}}}

\renewcommand{\Im}{\mathrm{Im}}
\renewcommand{\Re}{\mathrm{Re}}
\newcommand{\qedwhite}{\hfill \ensuremath{\Box}}

\newtheorem{remark}{Remark}[section]
\newtheorem{theorem}{Theorem}[section]
\newtheorem{lemma}{Lemma}[section]

\usepackage{xcolor,todonotes}
\author{M.~Schneider, J.~Lang, R.~Weiner}
\title{Super-Convergent Implicit-Explicit Peer Methods with Variable Step Sizes}
\author{
Moritz Schneider\\
{\small \it Technische Universit\"at Darmstadt, Dolivostra{\ss}e 15, 64293 Darmstadt, Germany}\\
{\small moschneider@mathematik.tu-darmstadt.de} \\ \\
Jens Lang\footnote{corresponding author}\\
{\small \it Technische Universit\"at Darmstadt, Dolivostra{\ss}e 15, 64293 Darmstadt, Germany}\\
{\small lang@mathematik.tu-darmstadt.de} \\ \\
R\"udiger Weiner \\
{\small \it Martin-Luther Universit\"at Halle-Wittenberg} \\
{\small \it Theodor-Lieser-Str. 5, 06120 Halle/Saale, Germany} \\
{\small ruediger.weiner@mathematik.uni-halle.de}}
\begin{document}
\maketitle

\begin{center}
\textit{Dedicated to the Memory of Willem Hundsdorfer (1954 - 2017)}\\[0.5cm]
\end{center}

\begin{abstract}
Dynamical systems with sub-processes evolving on many different time
scales are ubiquitous in applications. Their efficient solution is greatly enhanced by
automatic time step variation.
This paper is concerned with the theory, construction and application of
IMEX-Peer methods that are super-convergent for variable step sizes
and A-stable in the implicit part. IMEX schemes combine the necessary stability
of implicit and low computational costs of explicit methods to
efficiently solve systems of ordinary differential
equations with both stiff and non-stiff parts included in the source term.
To construct super-convergent IMEX-Peer methods which keep their higher
order for variable step sizes and exhibit favourable linear stability
properties, we derive necessary and sufficient conditions on
the nodes and coefficient matrices and apply an extrapolation approach based
on already computed stage values. New super-convergent IMEX-Peer methods of
order $s+1$ for $s=2,3,4$ stages are given as result of additional order
conditions which maintain the super-convergence property independent of step
size changes. Numerical experiments and a comparison to other super-convergent
IMEX-Peer methods show the potential of the new methods when applied with
local error control.

\end{abstract}

\noindent {\bf Keywords}: implicit-explicit (IMEX) Peer methods; super-convergence;
extrapolation; A-stability; variable step size; local error control

\section{Introduction}
Many dynamical processes in engineering, physics, chemistry and other areas
are modelled by large systems of ordinary differential equations (ODEs) of the form
\begin{equation}
u'(t) = F_0(u(t)) + F_1(u(t)),
\end{equation}
where $F_0:\R^m\rightarrow\R^m$ represents the non-stiff or mildly stiff part and
$F_1:\R^m\rightarrow\R^m$ gives the stiff part of the equation. Such problems
often result from semi-discretized systems of partial differential equations with
diffusion, advection and reaction terms. Instead of applying a single explicit or
implicit method, an often more appropriate and efficient approach is to use the
decomposition of the right-hand side by treating only the $F_1$ contribution in an implicit fashion.
Thus, favourable stability properties of implicit schemes and the advantage of lower
costs for explicit schemes are combined to enhance the overall computational efficiency.
Since dynamical systems typically have sub-processes evolving on many different time scales,
a good ODE integrator should come with some adaptive error control, making
frequent step size changes over its own progress. In smooth regions, a few large steps
should speed up the integration, whereas many small steps should be applied in
non-smooth terrains. The resulting gains in efficiency can be up to factors of
hundreds or more.

IMEX-Peer methods with variable step sizes have been successfully applied
by Soleimani, Knoth, and Weiner \cite{SoleimaniKnothWeiner2017}
to fast-wave-slow-wave problems arising in weather prediction.
The super-convergent IMEX-Peer methods recently developed by Soleimani and Weiner
\cite{SoleimaniWeiner2017a,SoleimaniWeiner2018} and Schneider, Lang,
and Hundsdorfer \cite{SchneiderLangHundsdorfer2018} can in principle be applied
with variable step sizes, but then they might lose their super-convergence property,
especially for serious step size changes. Super-convergent explicit Peer methods
for variable step sizes have first been constructed by Weiner, Schmitt, Podhaisky, and
Jebens \cite{WeinerSchmittPodhaiskyJebens2009}, exploiting special matrix structures.
Another approach to construct such methods is the use of extrapolation as proposed by
Schneider, Lang, and Hundsdorfer \cite{SchneiderLangHundsdorfer2018}. This idea goes back to Crouzeix
\cite{Crouzeix1980} and was also used by Cardone, Jackiewicz, Sandu, and Zhang
\cite{CardoneJackiewiczSanduZhang2014b,CardoneJackiewiczSanduZhang2014a} and later on
by Bra\'{s}, Izzo, and Jackiewicz \cite{BrasIzzoJackiewicz2017} to construct implicit-explicit
general linear and Runge-Kutta methods. The procedure can be easily extended to variable step
sizes for IMEX-Peer methods.

In this paper, we use the extrapolation approach to construct
new super-convergent IMEX-Peer methods that keep their higher order for variable step sizes
and exhibit favourable linear stability properties, including A-stability of the implicit
part. Additional order conditions on the nodes and coefficient matrices which maintain the
super-convergence property independent of step size changes are derived for
implicit, explicit and IMEX-Peer methods. We give formulas for new super-convergent
IMEX-Peer methods of order $s+1$ for $s=2,3,4$ stages. Stability regions are computed and compared
to those of super-convergent IMEX-Peer methods for constant step
sizes from Schneider, Lang, and Hundsdorfer \cite{SchneiderLangHundsdorfer2018}. Eventually,
numerical results are presented for a Prothero-Robinson problem,
the van der Pol oscillator, a one-dimensional Burgers equation with stiff diffusion
and a one-dimensional advection-reaction problem with stiff reactions.

\section{Implicit-Explicit Peer Methods with Variable Step Sizes}
\subsection{Super-convergent implicit Peer methods with variable
step sizes}\label{sec:impl-peer}
We apply the so-called Peer methods introduced by Schmitt, Weiner
and co-workers \cite{SchmittWeiner2004,SchmittWeiner2017,SoleimaniWeiner2017a} to
solve initial value problems in the vector space $\V=\R^m, m\ge 1$,
\begin{equation}
u'(t) = F(u(t)),\quad u(0)=u_0\in\V\,.
\end{equation}
The general form of an $s$-stage implicit Peer method with variable
step sizes $\st_n$ is
\begin{equation}
\label{implpeer1}
w_n = (P_n\otimes I)w_{n-1} + \st_n (Q_n\otimes I)F(w_{n-1}) +
\st_n (R_n\otimes I)F(w_n)
\end{equation}
with the $m\times m$ identity matrix $I$, the $s\times s$ coefficient matrices
$P_n=(p_{ij}(\sigma_n))$, $Q_n=(q_{ij}(\sigma_n))$, $R_n=(r_{ij}(\sigma_n))$, which
depend on the step size ratio $\sigma_n:=\st_n/\st_{n-1}$, and approximations
\begin{equation}
w_n = [w_{n,1}^T,\ldots,w_{n,s}^T]^T\in\V^s,\quad w_{n,i}\approx u(t_n+c_i\st_n).
\end{equation}
Here, $\V^s=\R^{ms}$, $t_n=\st_0+\ldots+\st_{n-1}$, $n\ge0$, and the nodes $c_i\in\R$ are
such that $c_i\ne c_j$ if $i\ne j$, and $c_s=1$.
Further, $F(w)=[F(w_i)]\in\V^s$ is the application of $F$ to all components
of $w\in\V^s$. The starting vector $w_0=[w_{0,i}]\in\V^s$ is supposed to be
given, or computed by a Runge-Kutta method, for example.

Peer methods belong to the class of general linear methods introduced
by Butcher \cite{Butcher2006}. All
approximations have the same order, which gives the name of the methods.
Here, we are interested in A-stable and super-convergent Peer methods with
order of convergence $p\!=\!s+1$ even for variable step sizes. For
constant step sizes, such methods have been recently constructed by Soleimani and
Weiner \cite{SoleimaniWeiner2017a} and Schneider, Lang and
Hundsdorfer \cite{SchneiderLangHundsdorfer2018}. In the following, for an
$s\times s$ matrix we will use the same symbol for its Kronecker product
with the identity matrix as a mapping from the space $\V^s$ to itself.
Then, (\ref{implpeer1}) simply reads
\begin{equation}
\label{implpeer}
w_n = P_nw_{n-1} + \st_n Q_nF(w_{n-1}) + \st_n R_nF(w_n)\,.
\end{equation}
In what follows, we discuss requirements and desirable
properties for the implicit Peer method (\ref{implpeer}).\\

\noindent\textbf{Accuracy.} Let $e=(1,\ldots,1)^T\in\R^s$. We assume
pre-consistency, i.e., $P_ne\!=\!e$, which means that for the trivial
equation $u'(t)=0$, we get solutions $w_{n,i}=1$ provided that
$w_{0,j}=1,\;j=1,\ldots,s$. The residual-type local errors result
from inserting exact solution values $w(t_n)=[u(t_n+c_i\st_n)]\in\V^s$
in the implicit scheme (\ref{implpeer}):
\begin{equation}
r_n = w(t_n) - P_n w(t_{n-1}) - \st_n Q_n w'(t_{n-1}) - \st_n R_n w'(t_n)\,.
\end{equation}
Let $c=(c_1,\ldots,c_s)^T$ with
point-wise powers $c^j=(c_1^j,\ldots,c_s^j)^T$. Then
Taylor expansion with the expressions
\begin{equation}
w_i(t_{n-1}) = u\left( t_n+\frac{c_i-1}{\sigma_n}\st_{n}\right),
\quad i=1,\ldots,s,
\end{equation}
gives
\begin{align}
w(t_n) = &\; e\otimes u(t_n) + \st_n c\otimes u'(t_n) +
  \frac{1}{2}\st_n^2 c^2\otimes u''(t_n) + \ldots\,\\
w(t_{n-1}) = &\; e\otimes u(t_n) + \frac{\st_n}{\sigma_n} (c-e)\otimes u'(t_n) +
  \frac{\st_n^2}{2\sigma_n^2} (c-e)^2\otimes u''(t_n) + \ldots,
\end{align}
from which we obtain
\begin{equation}
\label{def:res-r}
r_n = \sum_{j\ge 1} \st_n^j d_{n,j}\otimes u^{(j)}(t_n)
\end{equation}
with
\begin{equation}
\label{def:defect-d}
d_{n,j} = \frac{1}{j!}\left( c^j - \frac{1}{\sigma_n^j}P_n(c-e)^j
- \frac{j}{\sigma_n^{j-1}}Q_n(c-e)^{j-1} - jR_nc^{j-1}\right)\,.
\end{equation}
A pre-consistent method is said to have stage order $q$ if
$d_{n,j}=0$ for all $\sigma_n$ and $j=1,2,\ldots,q$. With the Vandermonde
matrices
\begin{equation}
\label{eq:V0V1}
V_0 = \big(c_i^{j-1}\big), \qquad
V_1 = \big((c_i-1)^{j-1}\big),\qquad i,j=1,\ldots,s,
\end{equation}
and $C=\text{diag}(c_1,c_2,\ldots,c_s)$,
$D=\text{diag}(1,2,\ldots,s)$, and
$S_n=\text{diag}(1,\sigma_n,\ldots,\sigma_n^{s-1})$,
the conditions for having stage order
$s$ for the implicit Peer method (\ref{implpeer}) for variable
step sizes are
\begin{equation}
CV_0 - \frac{1}{\sigma_n}P_n(C - I)V_1S_n^{-1} -
Q_nV_1 D S_n^{-1} - R_nV_0 D = 0\,.
\end{equation}
Since $V_1$ and $D$ are regular, we have the relation
\begin{equation}
\label{def:Q}
Q_n = \left( (CV_0 - R_nV_0 D)S_n - \frac{1}{\sigma_n}P_n(C - I)V_1
\right)(V_1D)^{-1}\,,
\end{equation}
showing that $Q_n$ is uniquely defined by the choice of $P_n$, $R_n$,
the node vector $c$, and the step size ratio $\sigma_n$. Moreover,
there is an easy way to achieve consistency for any choice of the step
sizes $\st_n$ by setting $P_n\equiv P$ and $R_n\equiv R$ with constant
matrices $P$ and $R$, and recomputing $Q_n$ from (\ref{def:Q}) in
each time step. In what follows, we will make use of this simplification
and consider implicit Peer methods with variable step sizes $\st_n$ of the form
\begin{equation}
\label{implpeervar}
w_n = Pw_{n-1} + \st_n Q_nF(w_{n-1}) + \st_n RF(w_n)
\end{equation}
with constant matrices $P$ and $R$, and $Q_n$ updated in each time step by
\begin{equation}
\label{def:Qvar}
Q_n = \left( (CV_0 - RV_0 D)S_n - \frac{1}{\sigma_n}P(C - I)V_1
\right)(V_1D)^{-1}\,.
\end{equation}
The matrix $R$ is taken to be lower triangular with constant diagonal
$r_{ii}\!=\!\gamma\!>\!0$, $i=1,\ldots,s$, giving singly diagonally
implicit methods.
\begin{remark}
Implicit Peer methods of the form (\ref{implpeervar}) that are consistent
of order $s$ for constant time steps, i.e., $\st_n=\st$ and $Q_n=Q$, can be
applied in a variable time-step environment without loss of their order of
consistency by updating (the original) $Q$ by $Q_n$ from (\ref{def:Qvar})
in each time step. We will use this modification in the numerical comparisons
for our recently developed methods in
\cite{SchneiderLangHundsdorfer2018,SoleimaniWeiner2018}.
\end{remark}

\noindent\textbf{Stability.}
Applying the implicit method (\ref{implpeervar}) to
Dahlquist's test equation $y'\!=\!\lambda y$ with $\lambda\in\C$,
gives the following recursion for the approximations $w_n$:
\begin{equation}
w_n = (I-z_{n}R)^{-1}(P+z_{n}Q_n)w_{n-1}=:M_{im}(z_{n},\sigma_{n})w_{n-1}
\end{equation}
with $z_{n}\!:=\!\sigma_n\sigma_{n-1}\cdots\sigma_1z_0,\;z_0=\st_0\lambda$.
Hence,
\begin{equation}
\label{eq:stab-matrix-func}
w_n=M_{im}(z_{n},\sigma_{n})M_{im}(z_{n-1},\sigma_{n-1})\cdots M_{im}(z_{1},\sigma_{1})w_{0}\,.
\end{equation}
The asymptotic behaviour of the matrix product is very difficult to
analyse, see e.g. the discussion in Jackiewicz, Podhaisky, and Weiner
\cite[Sect. 1]{JackiewiczPodhaiskyWeiner2004}. Here, we consider
methods that are zero-stable for arbitrary step sizes and A-stable
for constant step sizes. Zero stability requires the constant matrix
$P\!=\!M_{im}(0,\sigma)$ to be power bounded to have stability
for the trivial equation $u'(t)\!=\!0$. We will derive methods
for which the spectral radius of $M_{im}(z,\sigma)$ satisfies $\rho(M_{im}(z,\sigma))\le 1$ for all
$\sigma\in [\sigma_{min},\sigma_{max}]$ with $0\le\sigma_{min}<1\le\sigma_{max}$
and all $z\in\C$ with $\Re (z)\le 0$. Since for constant step sizes, $M_{im}(\infty,1)\!=\!-R^{-1}Q(1)$
with $Q(1)\neq 0$, A-stability does not imply L-stability. To guarantee good damping properties for very
stiff problems, we will aim at having a small spectral radius of $R^{-1}Q(\sigma)$ for
$\sigma\in [\sigma_{min},\sigma_{max}]$.
Although we cannot prove boundedness of the matrix product in (\ref{eq:stab-matrix-func}) for
$n\rightarrow \infty$ and variable step sizes,
the methods derived along the design principles described above performed always
stable in our numerical applications for various step size patterns.\\

\noindent\textbf{Super-convergence.}
Applying the convergence theory from multistep methods, stage order $q\!=\!s$ and
zero stability yield convergence of order $p\!=\!s$ for variable step sizes with
$0\le\sigma_{min}<\sigma_n<\sigma_{max}$ and $\st_n\le\st_{max}:=\max_{i=0,\ldots,n}\st_i$
demonstrated, e.g., in
\cite{BeckWeinerPodhaiskySchmitt2012,PodhaiskyWeinerSchmitt2005}.
Here, we are interested in using the degrees of freedom provided by
the free parameters in $P$, $R$, and $c$ to have convergence
of order $p\!=\!s+1$ without raising the stage order further.
This is discussed under the heading super-convergence in the book of
Strehmel, Weiner and Podhaisky \cite[Sect.\,5.3]{StrehmelWeinerPodhaisky2012}
for non-stiff problems. Similar results for stiff systems
were obtained by Hundsdorfer~\cite{Hundsdorfer1994}. We follow the
approach recently developed in Schneider, Lang, and
Hundsdorfer~\cite{SchneiderLangHundsdorfer2018}
for having an extra order of convergence for Peer methods with constant
step sizes to later discuss the property of super-convergence for
IMEX-Peer methods based on extrapolation for variable step sizes.

Let $\varepsilon_n\!=\!w(t_n)-w_n$ be the global error. Under the
standard stability assumption, where products of the transfer matrices
are bounded in norm by a fixed constant $K$ (see, e.g., Theorem 2 in
\cite{SoleimaniWeiner2017a}), we get the estimate
$\|\varepsilon_n\|\le K(\|\varepsilon_0\|+\|r_1\|+\ldots+\|r_n\|)$.
Together with stage order $s$, this gives the standard convergence result
\begin{align}
\label{implpeer-globerr}
\|\varepsilon_n\| \le &\; K\|\varepsilon_0\|+
K \left( \st_1^{s+1}\|d_{1,s+1}\|_\infty + \ldots +
\st_n^{s+1}\|d_{n,s+1}\|_\infty \right) \,\times\nonumber \\[1mm]
&\; \times\;\max_{0\le t\le t_n}\|u^{(s+1)}(t)\| +
\Oh\left(\st_{max}^{s+1}\right)\,.
\end{align}
Then we have the following
\begin{theorem}
\label{Th:super-implpeer}
Assume the implicit Peer method (\ref{implpeervar}) has stage order $s$
and estimate (\ref{implpeer-globerr}) holds true for the global error
with $\|\varepsilon_0\|=\Oh\left(\st_0^{s}\right)$.
Then the method is convergent of order $p\!=\!s$, i.e., the global
error satisfies $\varepsilon_n=\Oh\left(\st_{max}^s \right)$. Furthermore,
if the initial values are of order $s+1$, $d_{i,s+1}\in\text{range}\,(I-P)$
and $\st_{i-1}=\left( 1+\Oh(\st_{max})\right)\st_i$ for all $i=1,\ldots,n$,
then the order of convergence is $p\!=\!s+1$.
\end{theorem}
\noindent {\it Proof}: The first statement follows directly from
(\ref{implpeer-globerr}) with the estimate
\[ \st_1^{s+1}\|d_{1,s+1}\|_\infty + \ldots +
\st_n^{s+1}\|d_{n,s+1}\|_\infty \le
(t_{n+1}-t_1) \st_{max}^{s} \max_{i=1,\ldots,n}\|d_{i,s+1}\|_\infty\,.\]
Suppose that $d_{i,s+1}=(I-P)v_i$ with $v_i\in\R^s$. Since $I-P$ has
an eigenvalue zero, $v_i$ is not uniquely determined. To fix $v_i$, we
choose the one with minimum Euclidean norm, i.e., $v_i=(I-P)^+d_{i,s+1}$
with $(I-P)^+$ being the Moore-Penrose inverse. Let now
\begin{equation}
\bar{w}(t_i) := w(t_i) - \st_i^{s+1}v_i\otimes u^{(s+1)}(t_i)\,.
\end{equation}
Insertion of these modified solution values in the scheme
(\ref{implpeervar}) will give modified local errors
\begin{equation}
\label{eq:ri_bar}
\begin{array}{rll}
\bar{r}_i &=& \bar{w}(t_i)-P\bar{w}(t_{i-1})-\st_i Q_iF(\bar{w}(t_{i-1}))
-\st_i R F(\bar{w}(t_i))\\[2mm]
&=& r_i - \st_i^{s+1}d_{i,s+1}\otimes u^{(s+1)}(t_i)
- T(v_{i-1},v_i) \otimes u^{(s+1)}(t_i) + \Oh(\st_i^{s+2})\,,
\end{array}
\end{equation}
where
\begin{equation}
\label{eq:def_op_T}
T(v_{i-1},v_i) = \st_{i}^{s+1} Pv_i - \st_{i-1}^{s+1} Pv_{i-1}\,.
\end{equation}
Next, we will show that $T(v_{i-1},v_i)=\Oh(\st_{max})\st_i^{s+1}$. From
the assumption on the step sizes, $\st_{i-1}=(1+\Oh(\st_{max}))\st_{i}$, we deduce
$\sigma_i^{-1}=1+\Oh(\st_{max})$, which yields
\begin{equation}
\sigma_{i}^{-j} - \sigma_{i-1}^{-j} = \Oh(\st_{max})\quad\text{for all}\quad j\ge 1\,.
\end{equation}
The definition of $d_{i,s+1}$ in (\ref{def:defect-d}) gives the polynomial representation
$d_{i,s+1}=\sum_{j=0,\ldots,s+1}a_j\sigma^{-j}$ with $\sigma$-independent
$a_j\in\R^s$ (see also (\ref{eq:only_sigma}) and (\ref{eq:vd_sigma}) for more details).
Hence, we have $d_{i,s+1}-d_{i,s}=\Oh(\st_{max})$. Using
$\st_{i-1}^{s+1}=(1+\Oh(\st_{max}))\st_{i}^{s+1}$, we conclude that
\begin{equation}
\begin{array}{rll}
T(v_{i-1},v_i) &\!\!\!=\!\!\!& \st_{i}^{s+1}P(v_{i}-v_{i-1}) + \Oh(\st_{max})\st_i^{s+1}\\[2mm]
&\!\!\!=\!\!\!& \st_{i}^{s+1}P(I-P)^+(d_{i,s+1}-d_{i-1,s+1}) + \Oh(\st_{max})\st_i^{s+1}
= \Oh(\st_{max})\st_i^{s+1}\,,
\end{array}
\end{equation}
which, due to (\ref{def:res-r}), reveals $\bar{r}_i\!=\!\Oh(\st_{max}\st_i^{s+1})$ in (\ref{eq:ri_bar}).
This yields, in the same way as above,
$\|\bar{\varepsilon}_n\|\!=\!\|\bar{w}(t_n)-w_n\|\!\le\! K\|\varepsilon_0\|+\Oh(\st_{max}^{s+1})$.
Since $\|\bar{\varepsilon}_n-\varepsilon_n\|\!\le\! \st_n^{s+1}\|v_n\|_\infty\|u^{(s+1)}(t_n)\|$
and $\|\varepsilon_0\|\!=\!\Oh(\st_0^{s+1})$, this shows convergence of order $s+1$
for the global errors $\varepsilon_n$. \qedwhite \\

\noindent Recall that the range of $I-P$ consists of the vectors that are orthogonal to
the null space of $I-P^T$. If the method is zero-stable, then this null space has dimension one.
Therefore, up to a constant there is a unique vector $v\in\R^s$ such that $(I-P^T)v=0$. Then we have
\begin{equation}
\label{super-cond}
d_{i,s+1}\in\text{range}\,(I-P)\quad \text{iff}\quad
v^Td_{i,s+1}\!=\!0 \quad\text{for all }i=1,\ldots,n\,.
\end{equation}
Since $d_{i,s+1}$ depends on $\sigma_i$, these equations have to be satisfied
for all $\sigma_i$. In the following, we will drop the index $i$ and examine $v^Td_{s+1}(\sigma)$
as a function of $\sigma$. From (\ref{def:defect-d}), we
find
\begin{equation}
\label{eq:only_sigma}
d_{s+1}(\sigma) = \frac{1}{(s+1)!}
\left( c^{s+1} - \frac{1}{\sigma^{s+1}}P(c-e)^{s+1}
-\frac{s+1}{\sigma^s}Q(\sigma)(c-e)^s-(s+1)Rc^s\right),
\end{equation}
where $Q(\sigma)$ is taken from (\ref{def:Qvar}) with $\sigma_n=\sigma$.
Replacing $Q$, using the definition of $S_n$, and separating all powers of
$\sigma$, we eventually get the polynomial representation
\begin{equation}
\label{eq:vd_sigma}
v^Td_{s+1}(\sigma) = h_0 + \sum_{j=1}^{s}\tilde{v}_{s+1-j}^T
\,\tilde{c}_{s+1-j}\sigma^{-j}  + h_{s+1} \sigma^{-(s+1)}
\end{equation}
with the $\sigma$-independent coefficients
\begin{align}
h_0 = &\; \frac{1}{(s+1)!}v^T\left( c^{s+1} - (s+1)Rc^s \right),\\[2mm]
\tilde{v}^T = &\; \frac{1}{s!}v^T(RV_0-CV_0D^{-1}),\\[2mm]
\tilde{c} = &\; V_1^{-1}(c-e)^s,
\label{eq:tilde_c}\\[2mm]
h_{s+1} =&\; \frac{1}{(s+1)!}v^T (C-I)V_1\tilde{D}V_1^{-1}(c-e)^s.
\label{eq:h_s+1}
\end{align}
Here, $\tilde{D}:=(s+1)D^{-1}-I$. Note that we have used the relation
$v^TP=v^T$ to eliminate $P$ in (\ref{eq:h_s+1}). With
$h_j:=\tilde{v}_{s+1-j}^T\,\tilde{c}_{s+1-j}$ for $j=1,\ldots,s$, condition
(\ref{super-cond}) can be fulfilled by adding the $s+2$ additional equations
$h_j\equiv 0$ to the consistency conditions in order to achieve super-convergence
for variable step sizes. The special structure of the coefficients $\tilde{c}_{s+1-j}$
allows the following statement.
\begin{lemma}
Assume $c_1,\ldots,c_{s-1}<1$ with $s\ge 2$, $c_i\ne c_j$ for $i\ne j$, and $c_s=1$. Then,
$\tilde{c}_1=0$ and $\tilde{c}_2,\ldots,\tilde{c}_s\ne 0$.
\end{lemma}
\noindent {\it Proof}: The conditions on $c_i$ guarantee the regularity
of the Vandermonde matrix $V_1$. Let
$x_i:=c_i-1$, $i=1\ldots,s$. Then, we have $x^s=(c-e)^s$ and $V_1=(x_i^{j-1})$,
$i,j=1,\ldots,s$.
From (\ref{eq:tilde_c}), we deduce $V_1\tilde{c}=x^s$. The choice $c_s=1$ yields $x_s=0$
and hence $\tilde{c}_1=0$ from the last equation. Further, assumption
$c_i<1$ gives $x_1,\ldots,x_{s-1}<0$. This allows division by $x_i$, resulting
in the linear equations
\begin{equation}
\label{eq:lgs_c_tilde}
\begin{pmatrix}
1 & x_1 & \cdots & x_1^{s-2}\\
\vdots & \vdots & \ddots & \vdots\\
1 & x_{s-1} & \cdots & x_{s-1}^{s-2}
\end{pmatrix}\,
\begin{pmatrix}
\tilde{c}_2 \\ \vdots \\ \tilde{c}_{s}
\end{pmatrix} =
\begin{pmatrix}
x_1^{s-1} \\ \vdots \\x_{s-1}^{s-1}
\end{pmatrix}.
\end{equation}
Now, let us consider the polynomial of order $s-1$,
\begin{equation}
p(x) = x^{s-1} - \sum_{k=0}^{s-2} \tilde{c}_{k+2}\,x^k.
\end{equation}
Then, $p(x_i)=0$ is the $i$-th row of system (\ref{eq:lgs_c_tilde}) and
hence $x_1,\ldots,x_{s-1}$ are the $s-1$ roots of $p$, i.e.,
$p(x)=(x-x_1)\cdots (x-x_{s-1})$. The theorem of Vieta
shows
\begin{equation}
-\tilde{c}_{s+1-j} = (-1)^j\kappa_j\quad\text{with}\quad
\kappa_j = \sum_{1\le i_1<i_2<\cdots <i_j\le s-1}x_{i_1}\cdots x_{i_j},
\;j=1,\ldots,s-1.
\end{equation}
Since $x_i<0$ for all $i=1,\ldots,s-1$, we observe that all products in the
sum have the same number of factors which have one and the same sign, i.e.,
the sums cannot vanish. More precisely, $\text{sgn}(\kappa_j)=(-1)^j$
and hence $\text{sgn}(\tilde{c}_{s+1-j})=-1$,
which proves the statement.\qedwhite\\

\noindent We are now ready to formulate additional simplified conditions for the
super-conver\-gence of implicit Peer methods when they are applied with
variable step sizes.
\begin{theorem}
\label{Th:super-implpeer-cond}
Assume the implicit Peer method (\ref{implpeervar}) has stage order $s$
and estimate (\ref{implpeer-globerr}) holds true for the global error
with $\|\varepsilon_0\|=\Oh\left(\st_0^{s+1}\right)$.
Let $\st_{i-1}=(1+\Oh(\st_{max}))\st_{i}$ for all $i=1,\ldots,n$.
Then the method is convergent of order $p\!=\!s+1$, i.e., the global
error satisfies $\varepsilon_n=\Oh\left(\st_{max}^{s+1} \right)$, if
for all $v\in\R^s$ with $(I-P^T)v=0$, the
following additional conditions are satisfied:
\begin{align}
v^T (C-I)V_1\tilde{D}V_1^{-1}(c-e)^s = &\;0,\\[1mm]
v^T \left( c^j - jRc^{j-1} \right) = &\; 0,\quad j=2,\ldots,s+1.
\end{align}
\end{theorem}
\noindent {\it Proof}: Condition $v^Td_{s+1}(\sigma)=0$ requires $h_j=0$ for
$j=0,\ldots,s+1$ in (\ref{eq:vd_sigma}). Observe that $h_s=0$ is always
satisfied since $c_s=1$ and hence $\tilde{c}_1=0$. The property $\tilde{c}_{j}\ne 0$ for $j=2,\ldots,s$,
leads to $\tilde{v}^T_j=0$, which is equivalent to $v^T(c^j-jRc^{j-1})=0$.
The remaining conditions follow directly from $h_0=h_{s+1}=0$.
\qedwhite

\subsection{Super-convergent explicit Peer methods for variable steps sizes}
Super-convergent explicit Peer methods for variable step sizes with a special
structure of the matrix $P$ have first been constructed by Weiner, Schmitt,
Podhaisky, and Jebens \cite{WeinerSchmittPodhaiskyJebens2009}. A convenient way
to construct such methods for more general $P$ is the use of extrapolation as
proposed by Schneider, Lang, and Hundsdorfer \cite{SchneiderLangHundsdorfer2018}.
This idea goes back to Crouzeix \cite{Crouzeix1980} and was also used by
Cardone, Jackiewicz, Sandu, and Zhang~\cite{CardoneJackiewiczSanduZhang2014a}
to construct implicit-explicit diagonally implicit multistage integration methods.
The procedure can be easily extended to variable step sizes.

Assume that all approximations $w_{n,j}$ obtained from method
(\ref{implpeervar}) have stage order $s$. Then, we can use $w_{n-1}$ and
most recent values $w_{n,j}$, $j=1,\ldots,i-1$, already available for
the computation in the $i$-th stage, to extrapolate $F(w_n)$ by
\begin{equation}
F(w_n) = E_{1,n}F(w_{n-1}) + E_{2,n}F(w_n) + \Oh\left(\tau_n^{s}\right),
\end{equation}
where $\tau_n=\max{(\st_{n-1},\st_n)}$ and the $s\times s$-matrices
$E_{1,n}$ and $E_{2,n}$ of extrapolation
coefficients depend on the step size ratio $\sigma_n$. Here, $E_{2,n}$ is
a strictly lower triangular matrix. Replacing $F(w_n)$
in (\ref{implpeervar}) gives the explicit method
\begin{equation}
\label{explpeervar}
w_n = Pw_{n-1} + \st_n (Q_n + RE_{1,n})F(w_{n-1}) +
\st_n RE_{2,n}F(w_n).
\end{equation}
Note that $RE_{2,n}$ is strictly lower triangular since $R$ is
lower triangular. We will discuss consistency and super-convergence
of this explicit method.\\

\noindent\textbf{Accuracy.} Taylor expansion with exact values $F(w(t_n))$
gives for the residual-type error vector
\begin{equation}
\begin{array}{rll}
\delta_n &=& F(w(t_n))-E_{1,n}F(w(t_{n-1}))-E_{2,n}F(w(t_n)) \\[2mm]
 &=& \displaystyle\sum_{j\ge 0}\frac{\st_n^j}{j!}\left( (I-E_{2,n})c^j
 - \frac{1}{\sigma_n^j}E_{1,n}(c-e)^j
\right)\otimes\frac{d^j}{dt^j}F(u(t_n))\,.
\end{array}
\end{equation}
Then, the residual-type local error of the explicit Peer method
(\ref{explpeervar}) reads
\begin{equation}
r_n = \sum_{j\ge 1} \st_{n}^j \left( d_{n,j} + Rl_{n,j-1} \right)
\otimes u^{(j)}(t_n)
\end{equation}
with
\begin{equation}
\label{eq:expol_ln}
l_{n,j} = \frac{1}{j!}\left( (I-E_{2,n})c^{j}-\frac{1}{\sigma_n^j}E_{1,n}(c-e)^j\right).
\end{equation}
We can achieve stage order $s$, if the
underlying implicit Peer method has stage order~$s$, i.e.,
$d_{n,j}=0$ for all $\sigma_n$ and $j=1,\ldots,s$, and if
we choose
\begin{equation}
\label{eq:expol_E1n}
E_{1,n} = (I-E_2)V_0S_nV_1^{-1}
\end{equation}
with a constant $s\times s$-matrix $E_2$ and $S_n=\text{diag}(1,\sigma_n,\ldots,\sigma_n^{s-1})$ as defined
above. This gives $l_{n,j}=0$ for all $\sigma_n$ and $j=0,\ldots,s-1$
and eventually $r_n=\Oh(\st_n^{s+1})$.\\

\noindent\textbf{Super-convergence.} Under standard stability assumptions
as for the implicit method, we derive the global error estimate for the
explicit Peer method defined in (\ref{explpeervar}),
\begin{align}
\label{explpeer-globerr}
\|\varepsilon_n\| \le &\; K\|\varepsilon_0\|+
K \left( \st_1^{s+1}\|d_{1,s+1}+Rl_{1,s}\|_\infty + \ldots +
\st_n^{s+1}\|d_{n,s+1}+Rl_{n,s}\|_\infty \right) \,\times\nonumber \\[1mm]
&\; \times\;\max_{0\le t\le t_n}\|u^{(s+1)}(t)\| +
\Oh\left(\st_{max}^{s+1}\right)\,.
\end{align}
Analogously, we have
\begin{theorem}
\label{Th:super-explpeer}
Assume the implicit Peer method (\ref{implpeervar}) has stage order $s$
and estimate (\ref{explpeer-globerr}) holds true for the global error
with $\|\varepsilon_0\|=\Oh\left(\st_0^{s}\right)$.
Then the explicit method (\ref{explpeervar}) is convergent of order
$p\!=\!s$, i.e., the global
error satisfies $\varepsilon_n=\Oh\left(\st_{max}^s \right)$. Furthermore,
if the initial values are of order $s+1$,
$(d_{i,s+1}+Rl_{i,s})\in\text{range}\,(I-P)$ and
$\st_{i-1}=(1+\Oh(\st_{max}))\st_{i}$ for all $i=1,\ldots,n$, then the
order of convergence is $p\!=\!s+1$.
\end{theorem}
\noindent {\it Proof}: Replacing $d_{i,s+1}$ by $d_{i,s+1}+Rl_{i,s}$
in the proof of Theorem~\ref{Th:super-implpeer} gives the desired result.\qedwhite\\

\noindent Thus, super-convergence for variable step sizes is achieved if for all
$i=1,\ldots,n,$ it holds
\begin{equation}
v^T(d_{i,s+1}+Rl_{i,s})=0\text{ with }v\in\R^s\text{ such that }(I-P^T)v=0.
\end{equation}
If the underlying implicit method is already super-convergent, the
conditions simplify to $v^TRl_{i,s}=0$. Next, we will study the $l_{i,s}$
as functions of $\sigma$ and derive sufficient conditions for order $s+1$.

From (\ref{eq:expol_ln}) and (\ref{eq:expol_E1n}), we get
\begin{equation}
l_s(\sigma) = \frac{1}{s!} (I-E_2)\left( c^s -\frac{1}{\sigma^s}
V_0S(\sigma)V_1^{-1}(c-e)^s\right).
\end{equation}
The investigation of the product $v^T(d_{s+1}(\sigma)+Rl_s(\sigma))$
yields the following
\begin{theorem}
\label{Th:super-explpeer-var}
Assume the explicit Peer method (\ref{explpeervar}) has stage order $s$
and estimate (\ref{explpeer-globerr}) holds true for the global error
with $\|\varepsilon_0\|=\Oh\left(\st_0^{s+1}\right)$.
Let $\st_{i-1}=(1+\Oh(\st_{max}))\st_{i}$ for all $i=1,\ldots,n$.
Then the method is convergent of order $p\!=\!s+1$, i.e., the global
error satisfies $\varepsilon_n=\Oh\left(\st_{max}^{s+1} \right)$, if
for all $v\in\R^s$ with $(I-P^T)v=0$, the
following additional conditions are satisfied:
\begin{align}
\label{th:expl-cond1}
v^T (C-I)V_1\tilde{D}V_1^{-1}(c-e)^s = &\;0,\\[1mm]
\label{th:expl-cond2}
v^T \left( c^j - jRE_2c^{j-1} \right) = &\; 0,\quad j=2,\ldots,s+1.
\end{align}
\end{theorem}
\noindent {\it Proof}: The proof follows the same way as demonstrated
in the proof of Theorem~\ref{Th:super-implpeer}. The coefficients of
$\sigma^{-s},\ldots,\sigma^{-1}$ are again expressed as
products $\tilde{v}_1\tilde{c}_1,\ldots,\tilde{v}_s\tilde{c}_s$
with
\[ \tilde{c}=V_1^{-1}(c-e)^s \quad\text{and}\quad
\tilde{v}^T= \frac{1}{s!}v^T(RE_2V_0 -CV_0D^{-1}).\]
Due to $\tilde{c}_j\ne 0$ for $j=2,\ldots,s$, we have $\tilde{v}^T_j=0$
and hence $v^T(c^j-jRE_2c^{j-1})=0$. The other condition remains unchanged.
\qedwhite\\

\noindent We would like to conclude with the following observation:
If we start with a
super-convergent implicit Peer method for variable
step sizes, i.e., the additional conditions in Theorem~\ref{Th:super-implpeer-cond}
are already fulfilled, then (\ref{th:expl-cond1}) disappears and (\ref{th:expl-cond2})
changes to $v^TR(E_2-I)c^{j-1}=0$.
This can be rewritten to $v^TR(E_2-I)CV_0=0$. Since $R(E_2-I)$ and $V_0$ are
regular matrices, $C$ must be singular to satisfy (\ref{th:expl-cond2})
for $v\ne 0$. That means, one of the
nodes $c_i$ must be zero, because we always assume $c_i\ne c_j$. We will discuss this point later.

\subsection{Super-convergent IMEX-Peer methods with variable step sizes}
We now apply the implicit and explicit methods (\ref{implpeervar}) and (\ref{explpeervar})
to systems of the form
\begin{equation}
\label{ode-imex}
u'(t) = F_0(u(t)) + F_1(u(t))\,,
\end{equation}
where $F_0$ will represent the non-stiff or mildly stiff part, and $F_1$ gives
the stiff part of the equation. The resulting IMEX scheme is
\begin{equation}
\label{imex-peer-var}
w_n = P w_{n-1} + \st_n \left( \hQ_nF_0(w_{n-1}) + \hR F_0(w_{n})
+ Q_nF_1(w_{n-1}) + R F_1(w_{n})\right),
\end{equation}
where $\hQ_n=Q_n+RE_{1,n}$, $\hR=RE_2$, and
extrapolation is used only on $F_0$. Combining the local consistency analysis
for both the explicit and implicit method, the residual-type local errors for
the IMEX-Peer methods have the form
\begin{equation}
\label{imexpeer-res-var}
r_n = \sum_{j\ge 1}\st_n^j \left( d_{n,j}\otimes u^{(j)}(t_n) +
R\,l_{n,j-1}\otimes\frac{d^j}{dt^j}F_0(u(t_n)) \right).
\end{equation}

\noindent\textbf{Super-convergence.}
In order to construct super-convergent IMEX-Peer methods of order
$s+1$ for variable step sizes, we have to impose consistency of
order $s$ and ensure that for all $v\in\R^s$ with $(I-P^T)v=0$ it holds
\begin{equation}
v^Td_{s+1}(\sigma)=0\quad\text{and}\quad v^TR\,l_s(\sigma)=0
\end{equation}
for all $\sigma$. We have the following
\begin{theorem}
\label{Th:imex-conv-var}
Let the $s$-stage implicit Peer method (\ref{implpeervar}) defined by
the coefficients $(c,P,Q_n,R)$, with $Q_n$ from (\ref{def:Qvar}), be zero-stable
and suppose its stage order is equal to $s$.
Let the initial values satisfy $w_{0,i}-u(t_0+c_i\st_0)=\Oh(\st_0^{s+1})$, $i=1,\ldots,s,$
and $\st_{i-1}=\left( 1+\Oh(\st_{max})\right)\st_i$, $i=1,\ldots,n$.
Then the IMEX-Peer method (\ref{imex-peer-var}) is convergent
of order $s+1$, i.e., the global error satisfies $\varepsilon_n=\Oh(\st_{max}^{s+1})$,
if for all $v\in\R^s$ with $(I-P^T)v=0$,
the following additional conditions are satisfied:
\begin{align}
\label{th:imex-cond1}
v^T (C-I)V_1\tilde{D}V_1^{-1}(c-e)^s = &\;0,\\[1mm]
\label{th:imex-cond2}
v^T \left( c^j - jRc^{j-1} \right) = &\; 0,\quad j=2,\ldots,s+1,\\[1mm]
\label{th:imex-cond3}
v^T R(E_2-I)c^{j-1} = &\; 0,\quad j=2,\ldots,s+1.
\end{align}
\end{theorem}
\noindent {\it Proof}: Suppose $d_{i,s+1}=(I-P)v_{d,i}$ and $R\,l_{i,s}=(I-P)v_{l,i}$
with $v_{d,i},v_{l,i}\in\R^s$. Again, we fix these vectors by setting
$v_{d,i}=(I-P)^+d_{i,s+1}$ and $v_{l,i}=(I-P)^+Rl_{i,s}$ with $(I-P)^+$ being the
Moore-Penrose inverse. Let now
\begin{equation}
\bar{w}(t_i) = w(t_i) - \st_i^{s+1}v_{d,i}\otimes u^{(s+1)}(t_i) -
\st_i^{s+1}v_{l,i}\otimes \frac{d^s}{dt^s}F_0(u(t_i))\,.
\end{equation}
Inserting these modified values in (\ref{imex-peer-var}) gives the modified
residual-type local errors
\begin{equation}
\begin{array}{rll}
\label{}
\bar{r}_i &=& \bar{w}(t_i) - P \bar{w}(t_{i-1}) - \st_i\hQ_i F_0(\bar{w}(t_{i-1}))
- \st_i\hR F_0(\bar{w}(t_i))\\[2mm]
 && - \st_i Q_iF_1(\bar{w}(t_{i-1})) - \st_i R F_1(\bar{w}(t_i))\,,
\end{array}
\end{equation}
which can be rearranged to
\begin{equation}
\begin{array}{rll}
\bar{r}_i &=& \bar{w}(t_i) - P \bar{w}(t_{i-1}) - \st_i Q_i F(\bar{w}(t_{i-1}))
- \st_i R F(\bar{w}(t_i))\\[2mm]
&& + \st_i R \left( F_0(\bar{w}(t_i)) - E_{1,i}F_0(\bar{w}(t_{i-1}))
- E_2F_0(\bar{w}(t_i)) \right)\,.
\end{array}
\end{equation}
Then, Taylor expansions yields
\begin{equation}
\begin{array}{rll}
\bar{r}_i &=&
\displaystyle r_i - \st_i^{s+1}d_{i,s+1}\otimes u^{(s+1)}(t_i)
- \st_i^{s+1}R\,l_{i,s}\otimes \frac{d^s}{dt^s}F_0(u(t_i)) \\[3mm]
&& \displaystyle + \, T(v_{d,i-1},v_{d,i})\otimes u^{(s+1)}(t_i) + T(v_{l,i-1},v_{l,i})
\otimes \frac{d^s}{dt^s}F_0(u(t_i)) +  \Oh(\st_{i}^{s+2})
\end{array}
\end{equation}
with $T(\cdot,\cdot)$ and $r_i$ as defined in (\ref{eq:def_op_T}) and (\ref{imexpeer-res-var}),
respectively. The same arguments as in the proof of Theorem \ref{Th:super-implpeer} show
$\bar{r}_i=\Oh(\st_{max}\st_i^{s+1})$ and eventually the convergence of order $s+1$ for
the global errors $\varepsilon_n=w(t_n)-w_n$.\qedwhite\\

\noindent The $2s+1$ additional conditions (\ref{th:imex-cond1})-(\ref{th:imex-cond3}) are
quite demanding. We have already mentioned the fact that (\ref{th:imex-cond3})
requests that one of the nodes $c_i$, $i\ne s$, must be zero. In this case,
the method delivers two vectors, $w_{n-1,s}$ and $w_{n,i}$
with a certain $i$, that approximate
$u(t_n)$. We note that the difference of these approximations is used in
the extrapolation process as an additional degree of freedom. The matrix
$E_{1,n}$ in (\ref{eq:expol_E1n}) is still well defined. However, it is
not always possible to construct such methods at all or with good stability
properties in particular. In many practical applications, it might be
sufficient that
the explicit method has the property of super-convergence for
variable step sizes and the implicit method is only super-convergent for
constant step sizes. We have constructed such methods as well. They have
to fulfill the following additional conditions for all $v\in\R^s$ with
$(I-P^T)v=0$ and for all $\sigma$:
\begin{equation}
v^Td_{s+1}(1)=0\quad\text{and}\quad v^T(d_{s+1}(\sigma)+Rl_s(\sigma))=0.
\end{equation}
Due to the second condition for $\sigma=1$, the first one can be replaced
by the often simpler requirement $v^TRl_s(1)=0$. Using Theorem~\ref{Th:super-explpeer-var}
and the definition of $Rl_s$, we
find the explicit relations
\begin{align}
\label{th:imex-cond1-varex}
v^T R(I-E_2) \left( c^s - V_0V_1^{-1}(c-e)^s \right) = &\;0,\\[1mm]
\label{th:imex-cond2-varex}
v^T (C-I)V_1\tilde{D}V_1^{-1}(c-e)^s = &\;0,\\[1mm]
\label{th:imex-cond3-varex}
v^T \left( c^j - jRE_2c^{j-1} \right) = &\; 0,\quad j=2,\ldots,s+1.
\end{align}
Compared to (\ref{th:imex-cond1})-(\ref{th:imex-cond3}), the number of
conditions has been significantly reduced. Moreover, since condition (\ref{th:imex-cond3})
disappeared, the restriction $c_i=0$ for a certain $i$ is
no longer necessary.

\subsection{Stability of super-convergent IMEX-Peer methods}
We consider the usual split scalar test equation
\begin{equation}
\label{imex-test-eq}
y'(t) = \lambda_0 y(t) + \lambda_1y(t),\quad t\ge 0,
\end{equation}
with complex parameters $\lambda_0$ and $\lambda_1$.
Applying an IMEX-Peer method (\ref{imex-peer-var}) to
(\ref{imex-test-eq}) gives the recursion
\begin{equation}
\label{imex-peer-lin}
w_{n} = \left( I-z^{(n)}_0\hR-z^{(n)}_1R \right)^{-1}
\left( P+z^{(n)}_0\hQ_n
+ z^{(n)}_1Q_n \right)w_{n-1} =: M_n(z^{(n)}_0,z^{(n)}_1)w_{n-1}
\end{equation}
with $z_i^{(n)}=\st_n\lambda_i$, $i=0,1$. As for the implicit method itself,
an analysis of matrix products formed by $M_1M_2\cdots M_n$ would be
far too complicated. Therefore, we restrict ourselves to constant step
sizes and require
\begin{equation}
\label{imex-peer-stabmat}
\rho (M(z_0,z_1)) \le 1
\end{equation}
with $z_i=\st\lambda_i$, $i=0,1$. Then,
the stability regions of the IMEX-Peer method applied with
constant step sizes are defined by the sets
\begin{equation}
\St_\alpha = \{z_0\in\C: (\ref{imex-peer-stabmat})\text{ holds for any }
z_1 \in\C \text{ with }|\Im (z_1)|\le -\tan(\alpha)\cdot\Re(z_1) \}
\end{equation}
in the left-half complex plane for $\alpha \in [0^\circ,90^\circ]$.
Further, we define the stability region of the corresponding explicit
method (with constant step sizes) as
\begin{equation}
\St_E = \{z_0\in\C: \rho (M(z_0,0)) \le 1\}
\end{equation}
with the stability matrix $M(z_0,0)=(I-z_0\hR)^{-1}(P+z_0\hQ)$.
Efficient numerical algorithms to compute $\St_\alpha$ and $\St_E$ are
extensively described
in \cite{CardoneJackiewiczSanduZhang2014a,LangHundsdorfer2017}.

Our goal is to construct IMEX-Peer methods
for which $\St_E$ is large and $\St_E\backslash \St_\alpha$ is as small as
possible for angles $\alpha$ that are close to $90^\circ$. We will construct super-convergent
IMEX-Peer methods with A-stable implicit part for constant step sizes, i.e.,
the stability region $\St_{90^\circ}$ is non-empty. Concerning variable step sizes,
we follow the design principles already stated in the
stability discussion in Section~\ref{sec:impl-peer}.

\subsection{Practical Issues}
\noindent\textbf{Starting procedure.} In order to execute the first step
of the IMEX-Peer method (\ref{imex-peer-var}), we have to choose $t_1$,
$\st_0$, $\st_1$, and need to approximate the
$s$ initial values $w_{0,i}\approx u(t_1-(1-c_i)\st_0)$. For this, we
apply a suitable integration method with continuous output, e.g. a
Runge-Kutta or BDF scheme, on the interval $[t_0,t_0+\tau]$ with $\tau>0$.
The accuracy of the continuous numerical solution $\tilde{w}(t)$ can be controlled by
standard step size control or by choosing $\tau$ sufficiently small.
Denoting the minimum and maximum component of
the node vector $c$ by $c_{min}$ and $c_{max}$, respectively, we require
\begin{equation}
t_1-(1-c_{min})\st_0 = t_0 \quad\text{and}\quad
t_1-(1-c_{max})\st_0 = t_0 + \tau.
\end{equation}
This linear system for $t_1$ and $\st_0$ has the unique solution
\begin{equation}
t_1 = t_0 +\frac{1-c_{min}}{c_{max}-c_{min}}\tau \quad\text{and}\quad
\st_0 = \frac{1}{c_{max}-c_{min}}\tau.
\end{equation}
The initial values are now taken from
\begin{equation}
w_{0,i} := \tilde{w}(t_1-\st_0+c_i\st_0) =
\tilde{w}\left( t_0 + \frac{c_i-c_{min}}{c_{max}-c_{min}}\tau\right),
\;i=1,\ldots,s.
\end{equation}
Note that $w_{0,i}=u_0$ for index $i$ with $c_i=c_{min}$. Eventually,
we set $\st_1=\st_0$.\\

\noindent\textbf{Step size selection.} We extend the approach proposed by
Soleimani, Knoth, and Weiner in \cite{SoleimaniKnothWeiner2017} to locally
approximate $\st_n^su^{(s)}(t_n)$, which mimics the leading error term of
an embedded solution of order $s-1$. Let $F=F_0+F_1$ and define
\begin{equation}
est := \st_n \sum_{i=1}^s \left( \alpha_i F(w_{n,i}) + \beta_i F(w_{n-1,i}) \right)
\end{equation}
with $\alpha$ and $\beta$ determined through
\begin{equation}
\label{eq:est_weight}
\alpha^T = \delta (s-1)!\,e_s^TV_0^{-1}\quad\text{and}\quad
 \beta^T = (1-\delta)\sigma_n^{s-1}(s-1)!\,e_s^TV_1^{-1},
\end{equation}
where $e_s^T=(0,\ldots,0,1)$ and $\delta\in [0,1]$ is chosen as a weighting factor.
Then Taylor expansion of the exact solution shows the desired property:
\begin{align}
&\; \st_n \sum_{i=1}^s \left( \alpha_i u'(t_n+c_i\st_n)
+ \beta_i u'\left( t_n+\frac{c_i-1}{\sigma_n}\st_n
\right)\right)\nonumber\\[1mm]
= &\; \st_n \left( (\alpha^Te) u'(t_n) + \ldots + \frac{\st_n^{s-1}}{(s-1)!}(\alpha^Tc^{s-1})
u^{(s)}(t_n)\right. \\[1mm]
  &\; + \left. (\beta^Te) u'(t_n) + \ldots + \frac{\st_n^{s-1}}{\sigma_n^{s-1}(s-1)!}(\beta^T(c-e)^{s-1})
u^{(s)}(t_n)\right) + \Oh(\st_n^{s+1})\nonumber\\[1mm]
= &\; \st_n^{s}u^{(s)}(t_n) + \Oh(\st_n^{s+1}).
\end{align}
In our numerical experiments, we have discovered that the use
of old function values, i.e., $\delta=0$ in (\ref{eq:est_weight}),
works quite reliable for stiff and very stiff problems. For
mildly stiff problems, the choice $\delta=1$ often leads to a slightly
better performance. For our examples in Section~\ref{sec:NumExamples},
we will present results for $\delta=0$.

The new step size is computed by
\begin{equation}
\st_{new} =
\min \left( 1.2,\max \left( 0.8,0.9\,err^{-1/s}
\right) \right) \st_n
\end{equation}
with the weighted relative maximum error
\begin{equation}
err = \max_{i=1,\ldots,m} \frac{|est_i|}{atol\,+\,rtol\,
(\delta |w_{n,s,i}|+(1-\delta)|w_{n-1,s,i}|)}\,.
\end{equation}
In order to reach the time end point $T$ with a step of averaged
normal length, we adjust after each step size $\st_{new}$ to
$\st_{new} = (T-t_n)/\lfloor (1+(T-t_n)/\st_{new})\rfloor$.

Given an overall tolerance $TOL$, the step is accepted and the computation
is continued with $\st_{n+1}=\st_{new}$, if $err\le TOL$. Otherwise, the
step is rejected and repeated with $\st_{n}=\st_{new}$.

\section{Construction of super-convergent IMEX-Peer methods with variable
step sizes}
\subsection{The case $s=2$}
First, we have a negative result. With $c_1=0$, $c_2=1$, and pre-consistency $Pe=e$,
the coefficient matrices are
\begin{equation}
c=
\begin{pmatrix}
0 \\ 1
\end{pmatrix},\quad
P=
\begin{pmatrix}
p_1 & 1-p_1 \\ p_2 & 1-p_2
\end{pmatrix},\quad
R=
\begin{pmatrix}
\gamma & 0 \\ r_{21} & \gamma
\end{pmatrix},\quad
E_2=
\begin{pmatrix}
0 & 0 \\
e_{21} & 0
\end{pmatrix}.
\end{equation}
The first condition (\ref{th:imex-cond1}) for super-convergence reads
$(-1/2,0)\,v=0$, which gives, up to scaling, $v=(0,1)^T$. Then,
(\ref{th:imex-cond2}) reduces to $1-2\gamma=1-3\gamma=0$, which is
not possible for any $\gamma$.

Next we try to find methods that satisfy
(\ref{th:imex-cond1-varex})-(\ref{th:imex-cond3-varex}) with $c_1\ne 0$.
There are indeed candidates with $c_1=2/3$, $p_2=0$, $e_{12}=3/(4\gamma)$, and
$r_{21}=3/4-2\gamma$. The remaining parameters $p_1$ and $\gamma$ are
chosen such that the implicit part is A-stable and the stability regions of
the IMEX-method are optimized. Good results are obtained for the following
method:
\begin{equation}
c=
\begin{pmatrix}
\frac{2}{3} \\[2mm] 1
\end{pmatrix},\quad
P=
\begin{pmatrix}
-\frac{19}{20} & \frac{39}{20} \\[2mm] 0 & 1
\end{pmatrix},\quad
R=
\begin{pmatrix}
\frac{17}{20} & 0 \\[2mm] -\frac{19}{20} & \frac{17}{20}
\end{pmatrix},\quad
E_2=
\begin{pmatrix}
0 & 0 \\[2mm]
\frac{15}{17} & 0
\end{pmatrix}.
\end{equation}
We will refer to this method as IMEX-Peer2sve.\\

\begin{table}[h!]
\centering
{\small
\begin{tabular}{|l|cccr|c|cc|}
\hline\rule{0mm}{5mm}\hspace{-0.1cm}
IMEX- & $|\St_{90^\circ}|$ & $x_{max}$ & $|\St_{0^\circ}|$ & $y_{max}$ &
$\rho(R^{-1}Q)$ & $c_{im}$ & $c_{ex}$\\[1mm]
\hline\rule{0mm}{5mm}\hspace{-0.1cm}
Peer2s & $2.15$ & $-1.41$ & $4.47$ & $1.21$ & $0.128$ &
$2.37\,10^{-1}$ & $3.23\,10^{-1}$\\[2mm]
Peer2sve & $6.68\,10^{-5}$ & $-5.68\,10^{-3}$ & $0.14$ & $0.36$ & $0.863$ &
$1.94\,10^{-1}$ & $2.83\,10^{-1}$\\[1mm]
\hline\rule{0mm}{5mm}\hspace{-0.1cm}
Peer3s & $2.67$ & $-1.58$ & $6.11$ & $1.69$ & $0.552$ &
$1.24\,10^{-1}$ & $1.68\,10^{-1}$\\[2mm]
Peer3sv & $0.11$ & $-0.25$ & $0.55$ & $0.43$ & $0.254$ &
$2.29\,10^{-1}$ & $1.43\,10^{-1}$\\[1mm]
\hline\rule{0mm}{5mm}\hspace{-0.1cm}
Peer4s & $1.07$ & $-1.45$ & $4.39$ & $1.00$ & $0.542$ &
$6.42\,10^{-2}$ & $1.17\,10^{-1}$\\[2mm]
Peer4sve & $1.66$ & $-1.68$ & $3.11$ & $0.92$ & $0.118$ &
$2.02\,10^{-2}$ & $3.37\,10^{-2}$\\[2mm]
Peer4sv & $1.34\,10^{-3}$ & $-4.05\,10^{-2}$ & $0.63$ & $0.67$ & $0.632$ &
$7.47\,10^{-2}$ & $6.75\,10^{-2}$\\[1mm]
\hline
\end{tabular}
}
\parbox{13cm}{
\caption{\small
Size of stability regions $\St_{90^\circ}$ and $\St_{0^\circ}$,
$x_{max}(\St_{90^\circ})$ at the negative real axis,
$y_{max}(\St_{0^\circ})$ at the positive imaginary axis, spectral radius of $R^{-1}Q$,
and error constants $c_{im}=|d_{s+1}|$ and $c_{ex}=|R\,l_s|$ for super-convergent IMEX-Peer
methods, including those from \cite{SchneiderLangHundsdorfer2018}.}
\label{tab-stab-imex-super}
}
\end{table}

\subsection{The cases $s=3$ and $s=4$}
In order to construct super-convergent methods for variable
step sizes, we have to satisfy conditions (\ref{th:imex-cond1})-(\ref{th:imex-cond3})
for all $v\in\R^s$ with $(I-P^T)v=0$ and one of the nodes $c_i$ being zero.
A surprisingly simple choice is $c_1=0$ and $v=e_1$, which yields the validity of
(\ref{th:imex-cond2}) and (\ref{th:imex-cond3}). Then, equation (\ref{th:imex-cond1})
yields one condition for the remaining nodes. We find $c_2=0.5$ for $s=3$ and
$c_3=(5c_2-1)/(10c_2-5)$ for $s=4$. Furthermore, the first row of $P$ is $e_1$,
which goes along with pre-consistency. The value of $c_3$ and the remaining coefficients
of $P$, $R$ and $E_2$ are chosen in such a way that the implicit Peer methods
are A-stable and the IMEX-Peer methods exhibit good stability properties
and moderate error constants. This has been done using the Matlab-routine \textit{fminsearch},
where we included the desired
properties in the objective function and used random start values for the remaining
degrees of freedom. Different combinations of weights in the objective function
have been employed to select promising candidates which were then tested in various problems.
We will refer to the methods finally selected as IMEX-Peer3sv and IMEX-Peer4sv.

We have also constructed a $4$-stage IMEX-Peer method, denoted by
IMEX-Peer4sve, with the property that
the explicit method is super-convergent for variable step sizes and the
implicit method is only super-convergent for constant step sizes. In this case,
conditions (\ref{th:imex-cond1-varex})-(\ref{th:imex-cond3-varex}) must be satisfied,
where $c_i$, $i=1,2,3$, are still free parameters. We set $v=e_s$, which gives
(\ref{th:imex-cond2-varex}) since then $v^T(C-I)=0$. The additional degrees of freedom in
the nodes allow us to achieve greater stability regions and smaller error constants
compared to IMEX-Peer4sv. The method found is optimally zero-stable, i.e.,
one eigenvalue of $P$ equals one (due to pre-consistency) and the others
are zero.

\begin{figure}[t]
\setlength{\unitlength}{1cm}
\centering
\includegraphics[width=0.31\textwidth]{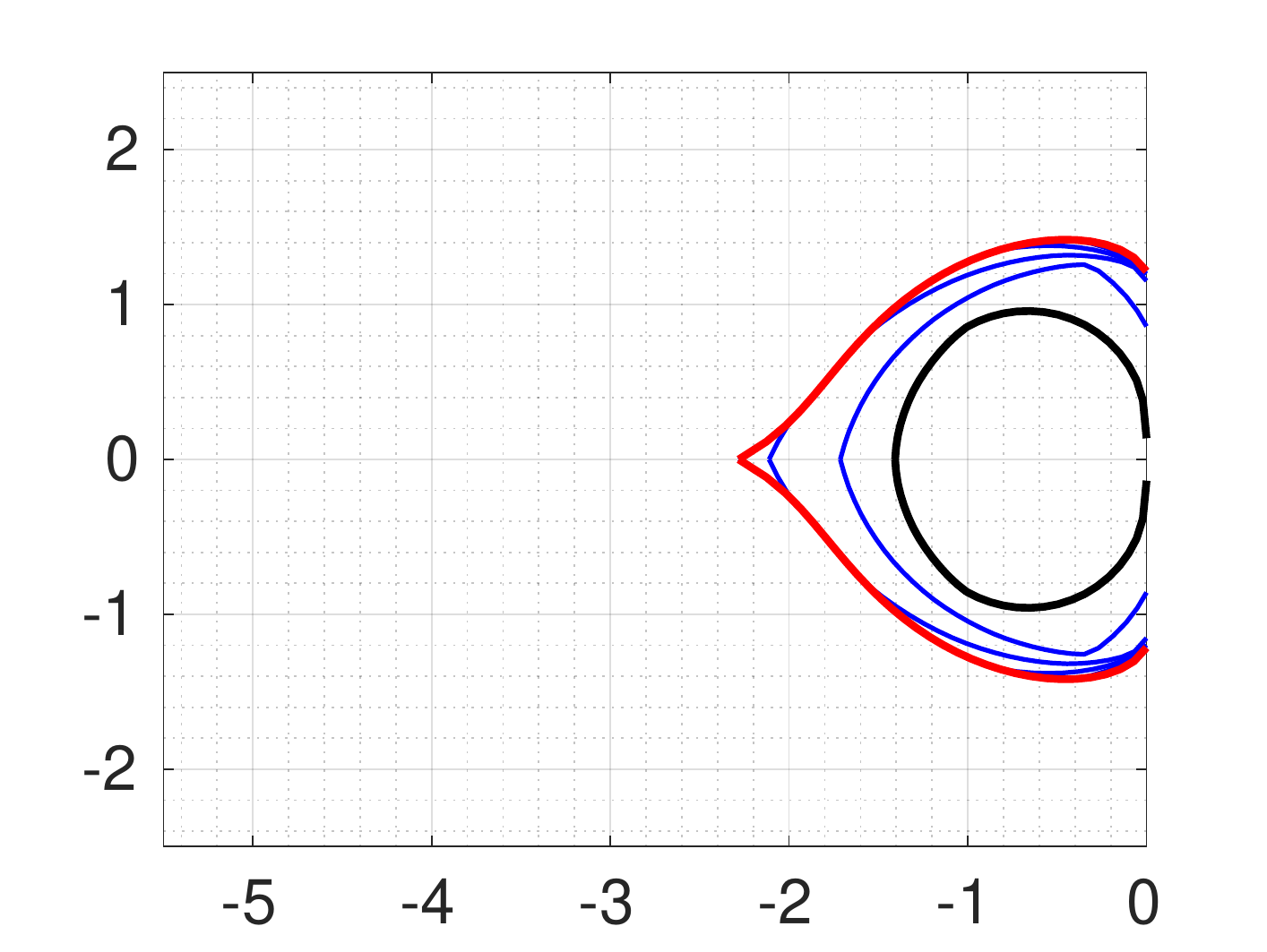}
\put(-3.6,2.7){\scriptsize IMEX-PEER2s}
\includegraphics[width=0.31\textwidth]{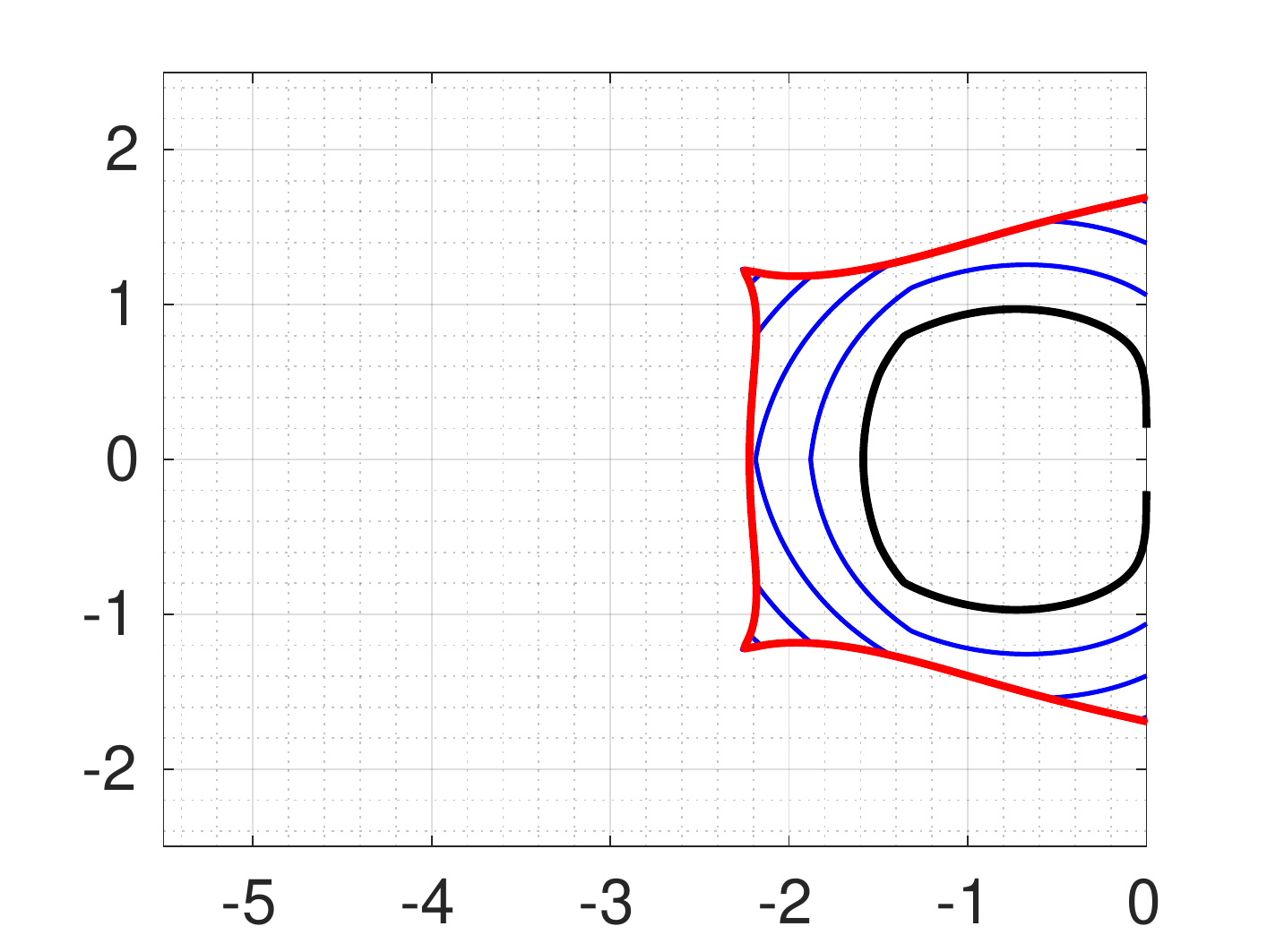}
\put(-3.6,2.7){\scriptsize IMEX-PEER3s}
\includegraphics[width=0.31\textwidth]{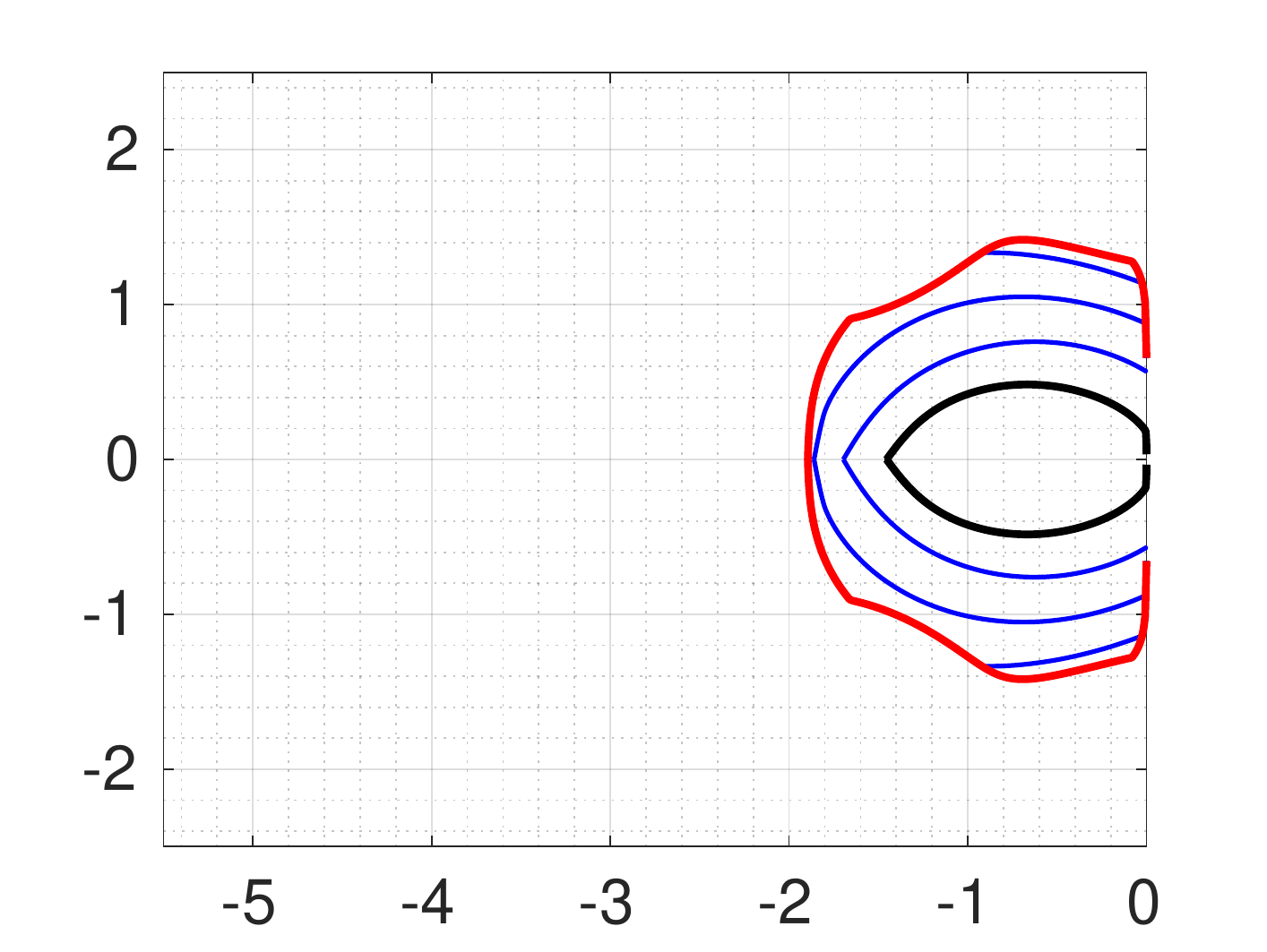}
\put(-3.6,2.7){\scriptsize IMEX-PEER4s}
\\
\includegraphics[width=0.31\textwidth]{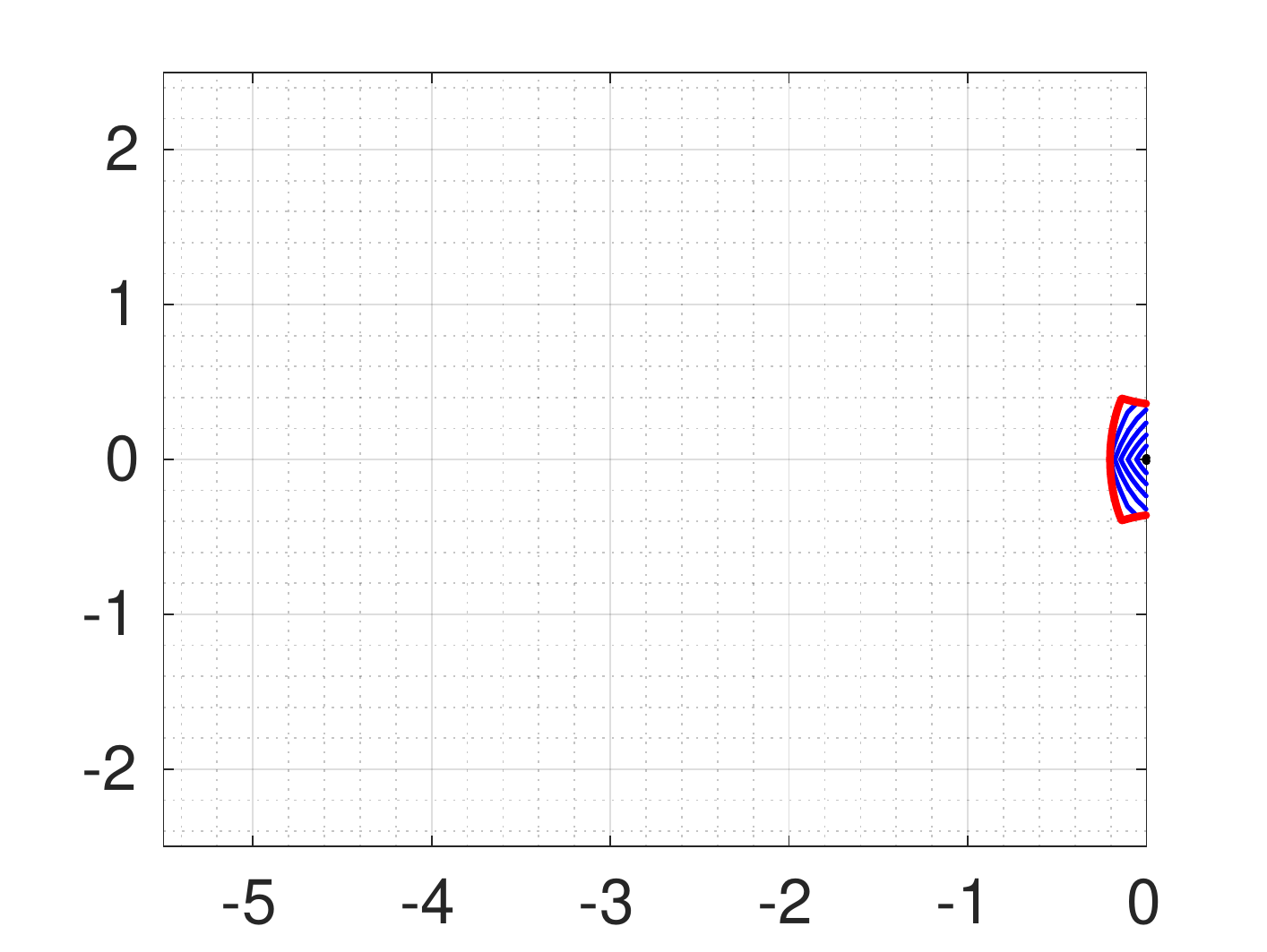}
\put(-3.6,2.7){\scriptsize IMEX-PEER2sve}
\includegraphics[width=0.31\textwidth]{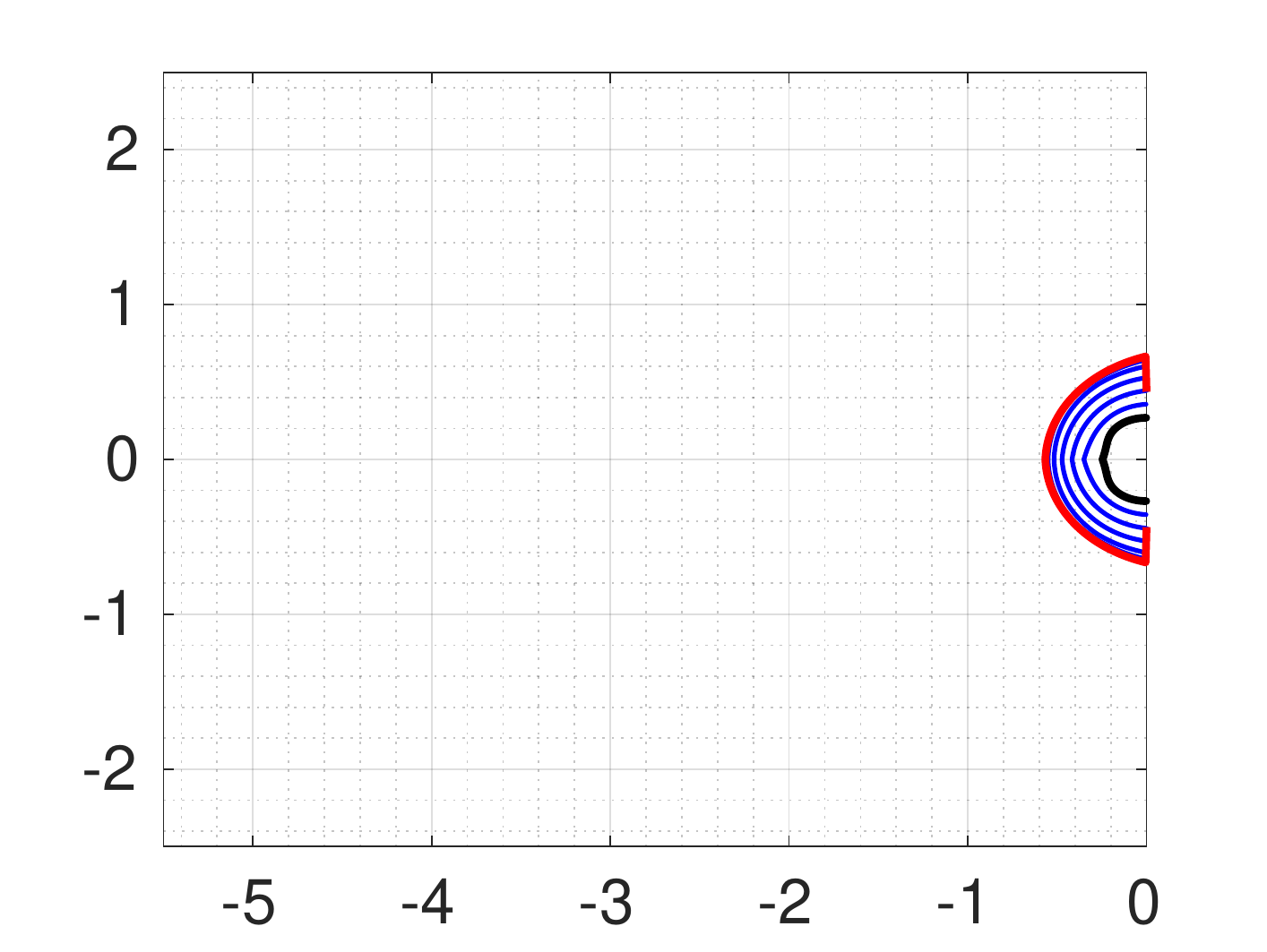}
\put(-3.6,2.7){\scriptsize IMEX-PEER3sv}
\includegraphics[width=0.31\textwidth]{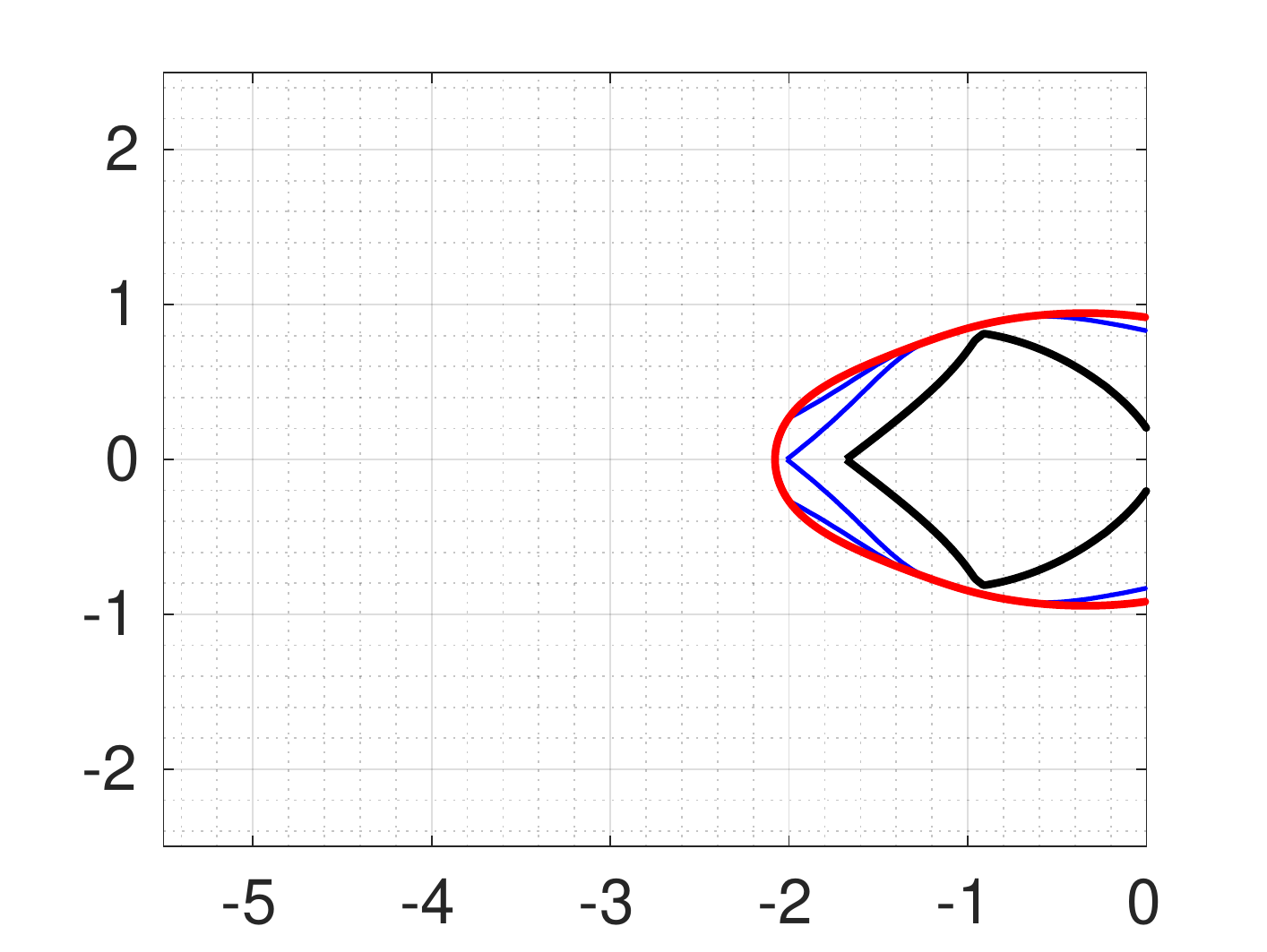}
\put(-3.6,2.7){\scriptsize IMEX-PEER4sve}
\\
\hspace{8.6cm}
\includegraphics[width=0.31\textwidth]{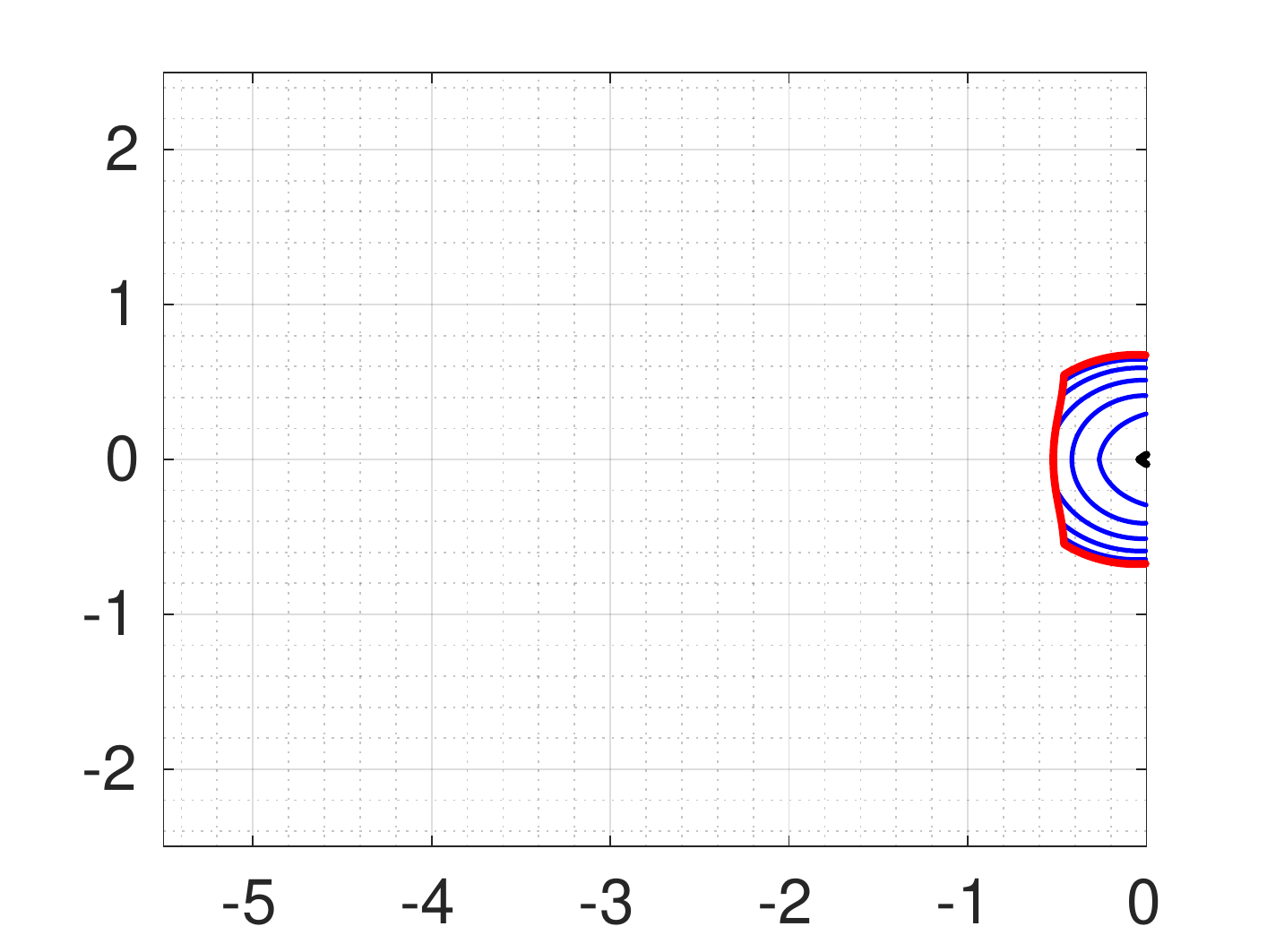}
\put(-3.6,2.7){\scriptsize IMEX-PEER4sv}
\\
\parbox{13cm}{
\caption{\small
Stability regions
$\St_{90^\circ}$ (black line), $\St_{\beta}$ for
$\beta=75^\circ,60^\circ,45^\circ,30^\circ,15^\circ$ (blue lines),
and $\St_{0^\circ}$ (red line) for super-convergent IMEX-Peer methods
with $s=2,3,4$ (left to right).}
\label{fig-stabreg-all-shape}
}
\end{figure}

The coefficients of all new methods for $c$, $P$, $R$, and $E_2$
are given in Table~\ref{tab:coeffs-34sv} and Table~\ref{tab:coeffs-4sve}.
Values for the stability regions as well as other constants are
collected in Table~\ref{tab-stab-imex-super}. More details on the stability
regions are shown in Figure~\ref{fig-stabreg-all-shape}. Obviously, the new property
of super-convergence for variable step sizes comes with significantly smaller
stability regions, except for IMEX-Peer4sve which even slightly improves
$S_{90^\circ}$ of IMEX-Peer4s.

\section{Numerical examples}
\label{sec:NumExamples}
We will present results for two ODE and two PDE problems. In order to
guarantee that errors of the initial values do not affect the computations,
unknown initial values as well as reference solutions $Y$ at the final time
are computed  by {\sc ode15s} from {\sc Matlab} with sufficiently high
tolerances. In the comparisons, the global errors are computed by
$err=\max_{i}|Y_i-\hat{Y}_i|/(1+|Y_i|)$, where $\hat{Y}$ is the numerical
approximation.

All calculations have been done with
Matlab-Version R2017a on a Latitude 7280 with an i5-7300U Intel processor
at 2.7 GHz.

\begin{table}[ht!]
\centering
{\footnotesize
\begin{tabular}{|l r c|c l r|}
\hline							
\multicolumn{6}{|l|}{IMEX-Peer3sv, $s=3$}\\
\hline
$c_1$	    &   $0.000000000000000$ 	&&&	$p_{11}$	& $1.000000000000000$ \\
$c_2$	    &   $0.500000000000000$ 	&&&	$p_{12}$	& $0.000000000000000$ \\
$c_3$	    &   $1.000000000000000$ 	&&&	$p_{13}$	& $0.000000000000000$ \\
$\gamma$    &   $0.690969692535085$     &&&	$p_{21}$	& $1.009534846612963$ \\
$r_{21}$    &   $0.351562922857064$     &&&	$p_{22}$	&$-0.000125189884283$ \\
$r_{31}$    &   $0.346024253990984$     &&&	$p_{23}$	&$-0.009409656728680$ \\
$r_{32}$    &   $0.328884660689640$     &&&	$p_{31}$	& $0.927244072163109$ \\
$e_{21}$    &   $1.454929231059714$     &&&	$p_{32}$	&$-0.000247968521087$ \\
$e_{31}$    &  $-6.099201725139450$     &&&	$p_{33}$ 	& $0.073003896357977$ \\		
$e_{32}$    &   $3.157746208382228$     &&& &   \\
\hline\hline
\multicolumn{6}{|l|}{IMEX-Peer4sv, $s=4$}\\
\hline
$c_1$		&   $0.000000000000000$     &&&     $p_{11}$ 	&  $1.000000000000000$ \\
$c_2$		&  $-1.598239239549169$     &&&	    $p_{12}$	&  $0.000000000000000$ \\
$c_3$		&   $0.523829503832339$     &&&	    $p_{13}$	&  $0.000000000000000$ \\
$c_4$		&   $1.000000000000000$     &&&	    $p_{14}$	&  $0.000000000000000$ \\
$\gamma$	&   $0.681884472048995$     &&&     $p_{21}$	&  $1.000204745561481$ \\
$r_{21}$	&   $1.292744499701930$     &&&     $p_{22}$	& $-0.000195233457439$ \\
$r_{31}$	&   $1.074957286644128$     &&&     $p_{23}$	& $-0.000009518220959$ \\
$r_{32}$	&  $-0.054028162784565$     &&&     $p_{24}$	&  $0.000000006116916$ \\
$r_{41}$	&   $4.064480810437903$     &&&     $p_{31}$	&  $1.169763235411655$ \\
$r_{42}$	&   $1.031994574173631$     &&&     $p_{32}$	& $-0.169740581681421$ \\
$r_{43}$	&  $-0.534558192336057$     &&&     $p_{33}$	& $-0.000025123517333$ \\
$e_{21}$	&  $-0.153830152235951$     &&&     $p_{34}$	&  $0.000002469787099$ \\
$e_{31}$	&   $0.065444441626366$     &&&     $p_{41}$	&  $1.915153835547942$ \\
$e_{32}$	&  $-0.976514386415223$     &&&     $p_{42}$	& $-0.244331567248295$ \\
$e_{41}$	&  $-0.234155732816782$     &&&     $p_{43}$	& $-0.671042624270695$ \\
$e_{42}$	&  $-2.535629358626096$     &&&     $p_{44}$	&  $0.000220355971049$ \\
$e_{43}$	&   $1.477107513945526$     &&&  & \\
\hline
\end{tabular}
}
\parbox{13cm}{
\caption{\small Coefficients of IMEX-Peer3sv and IMEX-Peer4sv
which are super-convergent for variable step sizes. Here, $E_2=(e_{ij})$.}
\label{tab:coeffs-34sv}
}
\end{table}

\begin{table}[ht!]
\centering
{\footnotesize
\begin{tabular}{|l r c|c l r|}
\hline
\multicolumn{6}{|l|}{IMEX-Peer4sve, $s=4$, optimally zero-stable}\\
\hline
$c_1$		&  $-0.868838855210029$     &&&     $p_{11}$ 	&  $0.000000000000000$ \\
$c_2$		&  $-0.253884413463736$     &&&	    $p_{12}$	&  $0.316402904545681$ \\
$c_3$		&   $0.754504864110948$     &&&	    $p_{13}$	&  $1.127642509582261$ \\
$c_4$		&   $1.000000000000000$     &&&	$    p_{14}$	& $-0.444045414127942$ \\
$\gamma$	&   $0.473861788489939$     &&&     $p_{21}$	&  $0.000000000000000$ \\
$r_{21}$	&   $0.732961380396538$     &&&     $p_{22}$	&  $0.000000000000000$ \\
$r_{31}$	&  $-2.472299983846101$     &&&     $p_{23}$	& $-0.017465269321373$ \\
$r_{32}$	&   $0.077358285702625$     &&&     $p_{24}$	&  $1.017465269321373$ \\
$r_{41}$	&  $-1.603925020256191$     &&&     $p_{31}$	&  $0.000000000000000$ \\
$r_{42}$	&  $-2.797576519478004$     &&&     $p_{32}$	&  $0.000000000000000$ \\
$r_{43}$	&  $-0.278164642408456$     &&&     $p_{33}$	&  $0.000000000000000$ \\
$e_{21}$	&  $-0.183287385063759$     &&&     $p_{34}$	&  $1.000000000000000$ \\
$e_{31}$	&   $5.974911797174020$     &&&     $p_{41}$	&  $0.000000000000000$ \\
$e_{32}$	&  $-2.556627399170977$     &&&     $p_{42}$	&  $0.000000000000000$ \\
$e_{41}$	&   $2.456065798975378$     &&&     $p_{43}$	&  $0.000000000000000$ \\
$e_{42}$	&  $-2.032396276261657$     &&&     $p_{44}$	&  $1.000000000000000$ \\
$e_{43}$	&   $1.255044479285407$     &&&  & \\
\hline
\end{tabular}
}
\parbox{13cm}{
\caption{\small Coefficients of IMEX-Peer4sve which is optimally zero-stable,
super-convergent for variable step sizes in the explicit part and for
constant step sizes in the implicit part. Here, $E_2=(e_{ij})$.}
\label{tab:coeffs-4sve}
}
\end{table}

\subsection{Prothero-Robinson Problem}
\begin{figure}[t!]
\setlength{\unitlength}{1cm}
\centering
\includegraphics[width=0.8\textwidth]{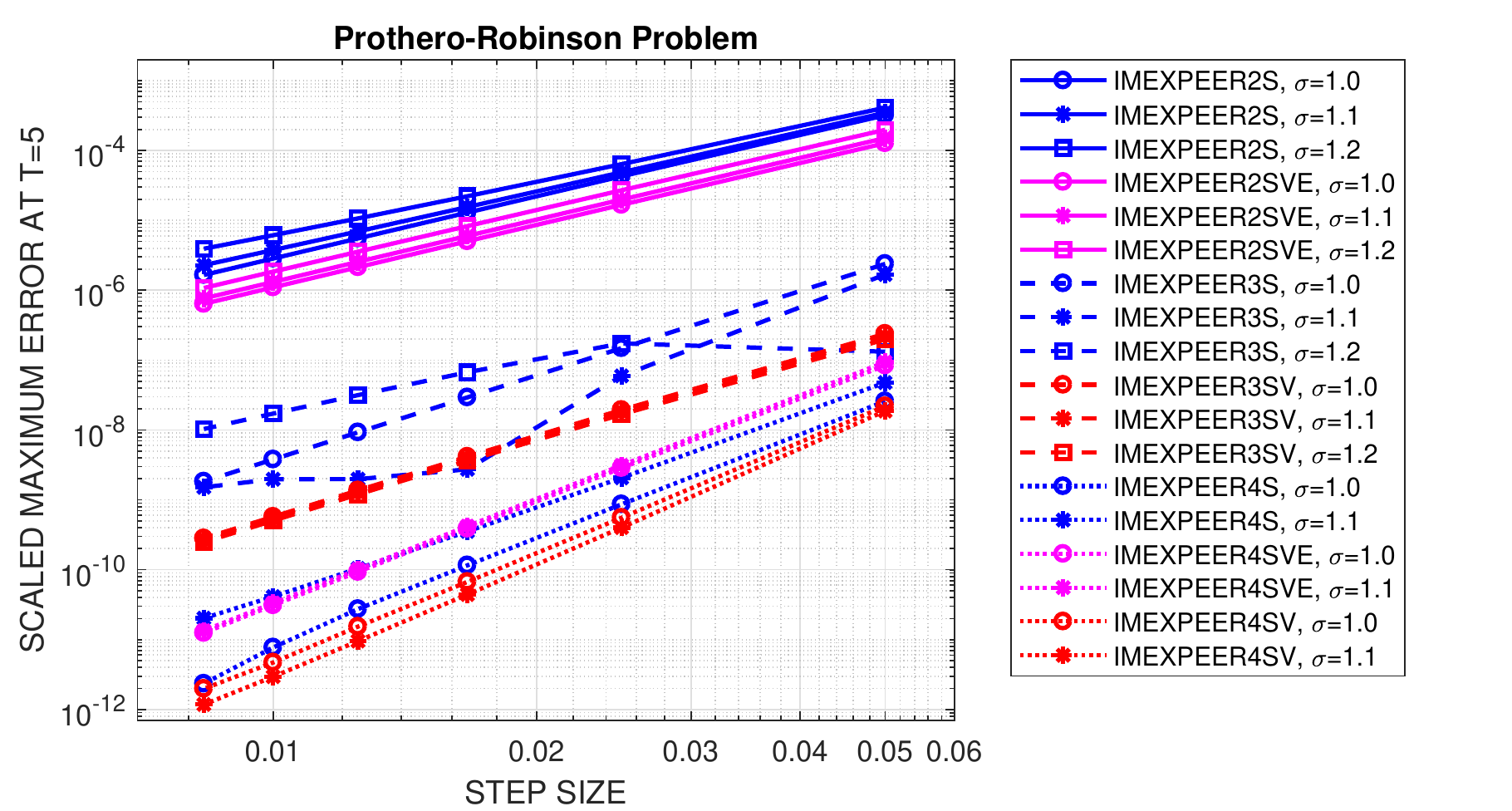}
\parbox{13cm}{
\caption{\small Prothero-Robinson Problem: Scaled maximum errors at $T=5$ vs.
time step sizes for $\sigma=1.0,1.1,1.2$. All new methods perform quite robust
with respect to changes of the step sizes.}
\label{fig:res-prothrob-sigma}
}
\end{figure}

In order to study the rate of convergence under stiffness and
changing step sizes, we consider the Prothero-Robinson type
equation used in \cite{SchneiderLangHundsdorfer2018,SoleimaniKnothWeiner2017},
\begin{align}
y' & =
\begin{pmatrix}
0 \\ y_1 + y_2 - \sin(t) \end{pmatrix} + \begin{pmatrix}
-10^6 (y_1 - \cos(t)) + 10^3 (y_2 -\sin(t)) - \sin(t) \\ 0
\end{pmatrix}\,,
\end{align}
where $t\in [0,5]$. The first term is treated explicitly and the second
implicitly. Initial values are taken from the analytic solution
$y(t)\!=\!(\cos(t),\sin(t))$. For constant step
sizes $\st=0.05/i$, $i=1,\ldots,6$, we consider the $\sigma$-dependent sequences
\begin{equation}
\st_i = \st_{i-1}\,\sigma^{(-1)^i},\quad i=2,\ldots,N
\end{equation}
with $\st_1=2\st/(1+\sigma)$ and $N=T/\st$. Results for $\sigma=1.0,1.1,1.2$
are shown in Figure~\ref{fig:res-prothrob-sigma}. Since the $4$-stage
methods become instable for $\sigma=1.2$, these results are omitted.
One can nicely see that all new methods keep their order of convergence
observed for constant step sizes and, therefore, perform quite robust
with respect to changing the step size. This is, of course, not the case
for the methods that are only super-convergent for constant
step sizes.

\subsection{Van der Pol Oscillator}
\begin{figure}[ht]
\setlength{\unitlength}{1cm}
\centering
\includegraphics[width=0.55\textwidth]{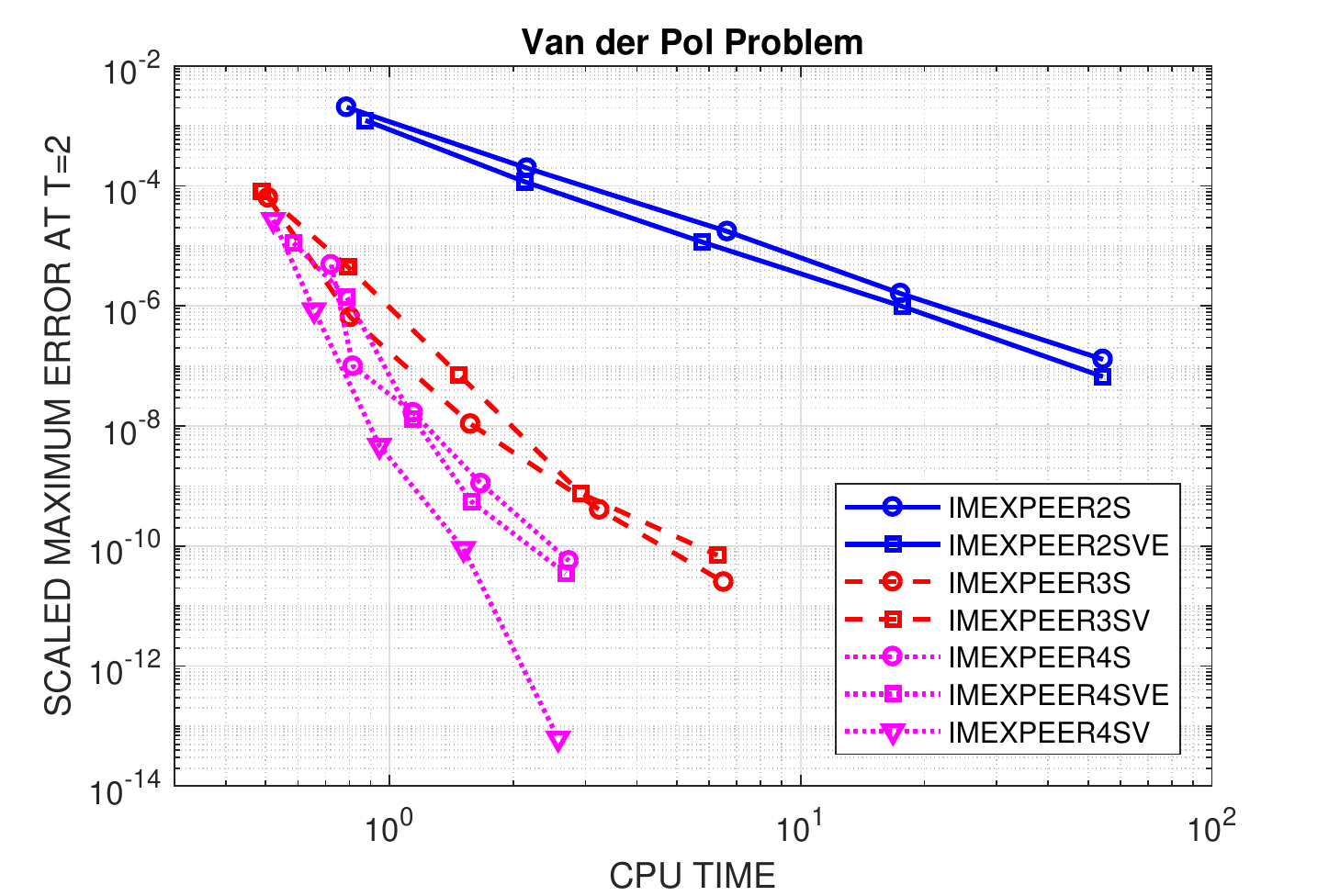}
\parbox{13cm}{
\caption{\small Van der Pol Oscillator: Scaled maximum errors at $T=2$ vs.
computing time.  For the $2$- and $3$-stage methods, the differences are
moderate. IMEX-Peer4sv shows a clear improvement over the other $4$-stage
methods.}
\label{fig:res-vanderpol}
}
\end{figure}

Next we consider the well known stiff van der Pol oscillator
\begin{align}
y' & =
\begin{pmatrix}
y_2 \\ 0
\end{pmatrix} +
\begin{pmatrix}
0 \\ 10^6\,((1-y_1^2)y_2-y_1)
\end{pmatrix}
\end{align}
with $y_1(0)=2$, $y_2(0)=0$, and $t\in [0,2]$. The first term is treated
explicitly and the second implicitly. This singularly perturbed problem
challenges any code and its efficient solution requires a step size
adaptation over several orders of magnitude, see e.g. \cite{HairerWanner1996}
and the discussions therein. The tolerances are $atol=rtol=10^{-3-i}$,
$i=0,1,\ldots,4$ and the calculations are started with initial step $\tau=atol$
for all methods. The results are shown and discussed in Figure~\ref{fig:res-vanderpol}.

\subsection{Burgers Problem}
The first PDE problem is taken from \cite{CalvoDeFrutos2001}, see also
\cite{SoleimaniWeiner2018} for further numerical results with super-convergent
IMEX-Peer methods. We consider
\begin{align}
\partial_t u = &\, 0.1\,\partial_{xx}u + u\partial_xu +\varphi(t,x),
\quad -1\le x\le 1,\;0\le t\le 2
\end{align}
with initial value $u(0,x)=\sin(\pi(x+1))$
and homogeneous Dirichlet boundary conditions. The
source term is defined through
\begin{align}
\varphi(t,x) = r(x)\sin(t),
\quad
r(x) = \left\{
\begin{array}{lrl}
0, & -1\!\!\!\! & \le x \le -1/3\\
3(x+1/3), & -1/3\!\!\!\! & \le x \le 0\\
3(2/3-x)/2, & 0\!\!\!\! & \le x \le 2/3\\
0, & 2/3\!\!\!\! & \le x \le 1.
\end{array}
\right.
\end{align}
The spatial discretization is done by finite differences with $\triangle x=1/2500$.
We treat the diffusion implicitly and all other terms explicitly.

\begin{figure}[t!]
\setlength{\unitlength}{1cm}
\centering
\includegraphics[width=0.48\textwidth]{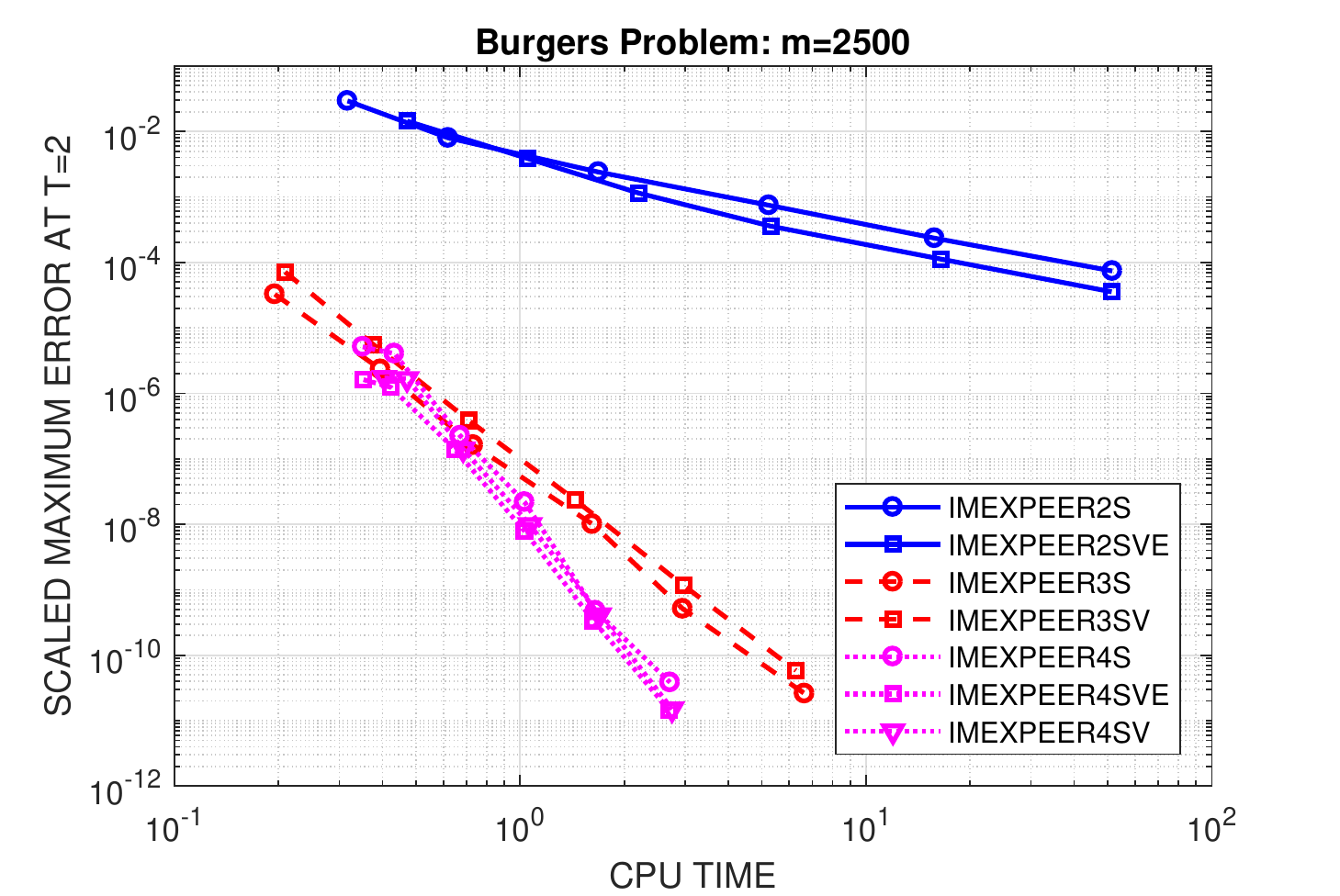}
\includegraphics[width=0.48\textwidth]{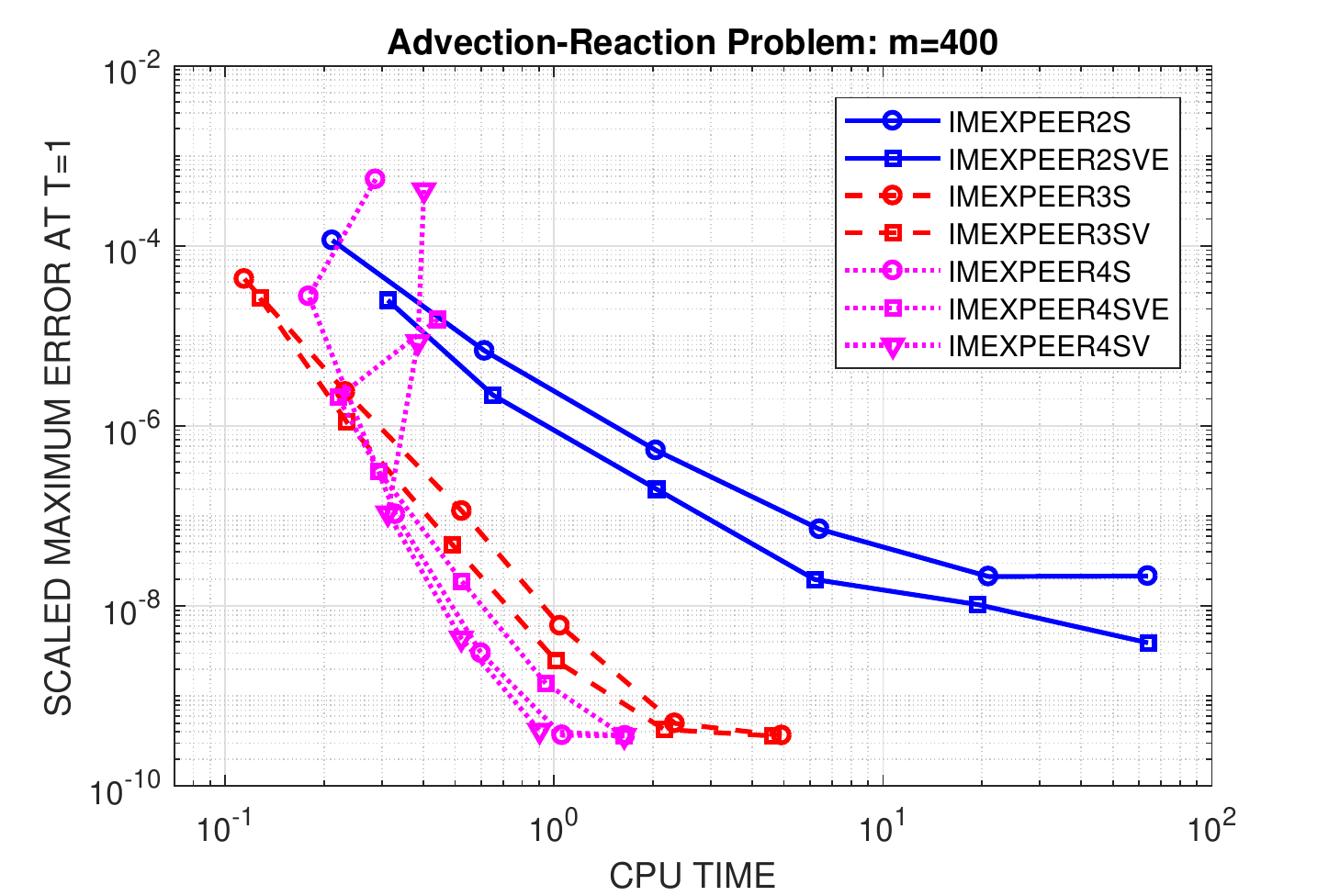}
\parbox{13cm}{
\caption{\small Burgers and Advection-Reaction Problem: Scaled maximum errors vs.
computing time. For the Burgers problem, no significant improvement can
be observed. In several cases, the better performance of the new methods for the
advection-reaction problem is obvious. All $4$-stage methods run for low tolerances
at their stability limit, which is related to $\st\approx 4\,10^{-4}$. The order
reduction of higher order methods for small time steps was already observed in
\cite{HundsdorferRuuth2007} and \cite{LangHundsdorfer2017}
as an inherent issue for very high-accuracy computations.}
\label{fig:res-burgers-advreac}
}
\end{figure}

We have used tolerances $atol=rtol=10^{-2-i}$, $i=0,1,\ldots,5$ and
initial step sizes $\tau=\sqrt{atol}$.
The results are plotted and discussed in Figure~\ref{fig:res-burgers-advreac}.

\subsection{Linear Advection-Reaction Problem}
A second PDE problem for an accuracy test is the linear advection-reaction system from
\cite{HundsdorferRuuth2007}. The equations are
\begin{eqnarray}
\partial_t u + \alpha_1\,\partial_x u &=& -k_1 u + k_2 v + s_1\,,\\
\partial_t v + \alpha_2\,\partial_x v &=& k_1 u - k_2 v + s_2
\end{eqnarray}
for $0<x<1$ and $0<t\le 1$, with parameters
\[
\alpha_1=1,\;\alpha_2=0,\;k_1=10^6,\;k_2=2k_1,\;s_1=0,\;s_2=1,
\]
and with the following initial and boundary conditions:
\[
u(x,0)=1+s_2x,\;v(x,0)=\frac{k_1}{k_2}u(x,0)+\frac{1}{k_2}s_2,\;
u(0,t)=1-\sin(12t)^4\,.
\]
Note that there are no boundary conditions for $v$ since $\alpha_2$
is set to be zero.

Fourth-order finite differences on a uniform mesh consisting of $m=400$
nodes are applied in the interior of the domain. At the boundary, we can
take third-order upwind biased finite differences, which here does not
affect an overall accuracy of four \cite{HundsdorferRuuth2007}
and gives rise to a spatial error of $1.5\,10^{-5}$.
In the IMEX setting, the reaction is treated implicitly and all other
terms explicitly.

We have used tolerances $atol=rtol=10^{-3-i}$, $i=0,1,\ldots,5$ and
an initial step size $\tau=10^{-3}$ for all runs.
The results are plotted and discussed in Figure~\ref{fig:res-burgers-advreac}.

\section{Conclusion}
We have developed a new class of $s$-stage super-convergent IMEX-Peer methods
with A-stable implicit part, which maintain their super-convergence order
of $s+1$ for variable step sizes. A-stability is important to solve problems with
function contributions that have large imaginary eigenvalues in the spectrum of
their Jacobian.
Applying the idea of extrapolation and studying the $\sigma$-dependent
coefficients in the local error representations, we first derived additional
conditions for implicit and explicit Peer methods, which are then combined to
state $2s+1$ corresponding conditions for IMEX-Peer methods. An interesting
theoretical result is that one of the nodes must be zero. Such methods exist
for $s\!>\!2$. We designed new methods for $s=3,4$. However, the new property
of super-convergence for variable step sizes reduces the scope for achieving
good stability properties, resulting in significantly smaller stability
regions compared to the super-convergent IMEX-Peer methods from
\cite{SchneiderLangHundsdorfer2018}. We also
constructed methods for $s=2,4$ having an explicit part that is super-convergent
for variable step sizes, whereas the implicit part is only super-convergent for
constant steps. In all cases, we employed the {\sc Matlab}-routine \textit{fminsearch}
with varying objective functions and starting values to find suitable methods
with stability regions as large as possible, good damping properties for very
stiff problems and small error constants.

We have implemented our newly designed methods with local error control based
on linear combinations of old function evaluations to approximate the
leading error term of an embedded solution of order $s\!-\!1$.
From our observations made for four numerical examples, we can draw the following
conclusions: (i) The new methods perform quite robust with respect to changing
the step size and, as expected, show their theoretical order at the same time.
(ii) For problems that demand a fast step size adaptation over several orders
of magnitudes, like the van der Pol oscillator, the new methods have the
potential to perform better. (iii) For problems that can be integrated with
moderate step size changes, like the Burgers problem, super-convergence for
constant step sizes is still sufficient to profit from the additional order and
possibly from the larger stability regions.

\section{Acknowledgement}
J.~Lang was supported by the
German Research Foundation within the collaborative research center
TRR154 ``Mathematical Modeling, Simulation and Optimisation Using
the Example of Gas Networks'' (DFG-SFB TRR154/2-2018, TP B01) and
the Graduate Schools Computational Engineering (DFG GSC233)
and Energy Science and Engineering (DFG GSC1070).

\bibliographystyle{plain}
\bibliography{bibimexpeersuper}

\end{document}